\documentclass[11pt,english,runninghead]{article}
\usepackage[T1]{fontenc}
\usepackage[latin9]{inputenc}
\usepackage[letterpaper]{geometry}
\usepackage{amsthm}
\usepackage{amsmath}
\usepackage{graphicx}
\usepackage{epstopdf}
\usepackage{setspace}
\usepackage{amssymb}
\usepackage{esint}

\makeatletter

\numberwithin{equation}{section}
  \theoremstyle{plain}
   
  \theoremstyle{remark}
  
  \theoremstyle{plain}

\usepackage{mathrsfs}

\makeatother

\usepackage{babel}

\makeatother

\usepackage{babel}

\usepackage[T1]{fontenc}

\usepackage{amsthm}
\usepackage{amsmath}
\usepackage{graphicx}
\usepackage{setspace}
\usepackage{amssymb}
\usepackage{esint}

\usepackage[final]{pdfpages}

\title{Dynamic  discrete models  for the   granular matter formation process}
\author{A. Khapalov and S. Lapin \\
Department of Mathematics \\
Pullman, WA 99163, USA \\
email: khapala@math.wsu.edu}
\date{}

\begin{document}

\maketitle

\bigskip

\begin{abstract}
In this paper we introduce 1-$D$ and 2-$D$ discrete  models for the  {\em dynamic} granular matter formation process in the form of a system of difference equations. This approach  allows us  to differentiate between  the influx   of the rolling layer coming  down from different directions to the corner points  of the standing layer. Such points are difficult to adequately describe by means of pde's and their straightforward numerical approximations, typically ``ignoring'' the system's behavior on the sets of zero measure. However,  these points are critical  for understanding the dynamic formation process when the standing layer is created by the moving front of the rolling matter or when the latter is filling a cavity and/or  stops rolling.  The existence of  distributed (infinite dimensional) limit solutions to our discrete models as the size of mesh  tends to zero is also discussed.  We illustrate our findings by numerical examples which use our model as a direct algorithm. 
\end{abstract}

{\bf 1. Introduction.}

{\bf 1.1.  Discussion of  prior work in the field.}  In the last  two decades there has been a substantial interest among physicists  and applied mathematicians to mathematical models, both dynamic and static, describing the process of formation of granular matter. The  motion patterns of granular materials are quite  different from the  behavior of  classical phases such as solid continua, gases or fluids, see the discussion in this respect, e.g.,  in de Gennes \cite{Gen1}-\cite{Gen2}, Hadeler and Kuttler \cite{hk}  and the references there in. 
Essentially, the formation of granular materials is viewed as the growth of so-called ``standing'' layer (whose height at point $ x$ at time $t$ we further denote by $ u(x,t) \geq 0 $) due to the influx  of a rolling layer (with the height $ v(x,t) \geq 0$), which rolls  down the slope of the standing layer under the gravity, see, e.g., Mehta et al \cite{Meh}, de Gennes \cite{Gen2} and the references there in. This process can be quite different for different  materials  (e.g., sand, snow, mud) due to the difference in their properties.  In this paper we focus on the case when the granular material is a dry sand or similar matter, in which case we do {\em not} deal with such phenomena like  mud slides or avalanches  over the  rock surface of the same fixed slope.

In Bouchaud et al \cite{Bch1}-\cite{Bch2} and  Boutreux et al \cite{Bou2}-\cite{Bou2}  the following 1-$D$ model was suggested for the {\em local} process of granular matter formation:
$$
v_t =  s v_x - \gamma (\alpha - \mid u_x \mid) v, 
\eqno(1.1)$$
$$
u_t = \gamma (\alpha - \mid u_x \mid)v \;\;\;\; {\rm in} \;\;  Q_T,
\eqno(1.2)$$
where  $s > 0 $ is to be the {\em constant} speed of the grains in the rolling layer when the standing layer slope is below its critical value $ \alpha$, and $ \gamma > 0 $ is a parameter. The choice of $ s>0$  means that the rolling motion is directed to the left. The  height of the standing layer at point $ x$ is to grow proportionally to the thickness of the rolling layer at this point. The height of the rolling layer grows proportionally to its slope minus the contribution to the standing layer.

Model  (1.1)-(1.2)  can be viewed as a ``minimal empiric''  model of the the granular pile formation process and it captured its essential features. However, it leaves open  a number  of questions and concerns. In particular, it assumes that  the speed of grains in the rolling layer does not depend on  the slope of the standing layer. This model  is also invalid when this slope  changes sign. No stochastic effects or inertia are taken into account. 
These and some other problems were  mentioned, e.g.,   in \cite{hk}.

Several subsequent papers were aimed at the attempts to address the above concerns. In \cite{Gen1} it was suggested to add a diffusion term $ D u_{xx} $ with small $ D$ to take into account possible random effects, and to add ``mathematical smoothness'' into model (1.1)-(1.2). In \cite{hk} the following  modifications of the above model were proposed: 
$$
v_t = \beta (v u_x)_x - \gamma (\alpha - \mid u_x \mid)v + f ,
\eqno(1.3)$$
$$
u_t = \gamma (\alpha - \mid u_x \mid)v,
\eqno(1.4)$$
$$
u(x,0) = u_o (x) \geq 0, \;\;  v(x,0) = v_o (x) \geq 0, \;\; \mid u_{ox} \mid \; \leq \alpha,
$$
where $ f $ describes a possible  external source of the granular matter. The first term on the right in (1.3) reflects the intention of the authors to make sure that the flow of the rolling layer follows the gradient of the standing layer. 

The above models  do not  discuss the question of conservation of granular matter.
Let us note that the above-cited models (both the original one and its modifications) were introduced in a rather  empiric fashion. In this paper we  intend to put more emphasis  on  the  physical and mathematical justification of modeling ideas  at hand.

In interesting papers Amadori and Shen \cite{Ama1}-\cite{Ama2} the existence of a global solution to system  (1.3)-(1.4) with $ f = 0$ was shown for small initial datum under the additional assumptions that the  slope of the standing layer does not change sign and  that it has an infinite length.  These assumptions convert  (1.3)-(1.4) into a qualitatively different  Cauchy problem for 2x2 hyperbolic  system of balance laws.  But they also make  the resulting model not quite realistic physically.

{\bf Other related publications.} In Prigozhin \cite{pr} (see also \cite{pr2} and the references therein) a variational approach to the granular matter formation problem was proposed in the form of a   model consisting of a single equation for the standing layer, while the height of the rolling layer is to be viewed as a ``selectable'' parameter. It was shown that this model is equivalent to an evolutionary quasi-variational inequality. The existence and uniqueness results  were established for the case of ``no steep'' slopes.

Several interesting explicit mathematical results were obtained in Cannarsa and Cardaliaguet et al {\cite{Can1}, Cannarsa et al\cite{Can2}-\cite{Can3}  for the {\em static version} of the aforementioned model in \cite{pr}.

{\bf 1.2. Motivation, main goals and layout of the paper.}  Our motivation for this paper is twofold.

Model (1.3)-(1.4) is well-recognized in the literature. However, recent  computational experiments   in Colombo et al  \cite{Col} revealed   certain  irregularities in its  behavior in some ``non-standard'' situations (such as when the standard layer forms  a cavity, see Example 2.1 below). It was noticed in \cite{Col} that the mathematics of the aforementioned model  equations allows the level of the granular matter to  exceed its initial maximum at the subsequent  moments of time.  This phenomenon -- labeled by  the authors of \cite{Col} as ``geysers'' -- seems to contradict  the main underlying physical idea that the rolling matter  should flow downward. Respectively,  \cite{Col} proposed   a  new very different model for the granular flow  in the form of  Cauchy problem for 3x3 system of balance laws.  
{\em In this paper} our goal is to try   to offer a different modification of the original flow model  (1.3)-(1.4). To this end, we investigate analytically which part  of equations  (1.3)-(1.4) is responsible  for generation of ``geysers'' and then modify this part (namely, by omitting the 2nd order term in (1.3)).

The 2nd  motivation point is that it seems that model (1.3)-(1.4) is not intended for the description of  the behavior of  the granular formation  process near the corners of the standing layer, that is, at the points where the function describing it is not differentiable.  
More precisely,  we  focus here on: (a)  the behavior of  granular matter near  the  points where  the  standing layer meets its foundation and its slope is critical (i.e., {\em how the slope propagates horizontally?}) and (b)
the  process of filling   cavities formed by  the standing layer with critical slopes.
{\em In this paper} we intend to overcome this challenge by proposing models in the form of {\em difference equations}, and show that,  as the mesh size tends to zero, the respective discrete solutions converge to some limit distributed function (along some subsequences, in general).  Our discrete models are linked to respective pde models under some regularity assumptions assumptions,  which, in particular,  exclude the aforementioned corner points. 

\bigskip
The paper is organized as follows.
In Section 2 we  prove  (by simple analytical means) that the  maximum value of the total height of the granular matter in the 1-$D$  model (1.3)-(1.4) may exceed its initial value  at some later time.  We show   that   this phenomenon is linked to the curvature of the surface of the standing layer  and does {\em not} depend on the thickness of the rolling layer, nor on the  thickness of the standing layer (see also Remark 2.1 below). The size of such emerging ``geysers'' can be substantial and cannot be viewed as just a minor effect of local diffusion, which was the reason in \cite{Gen1} to introduce the 2-nd order term $ \beta v u_{xx } $ into the 1-st equation of model (1.3)-(1.4). 
We show that, if we  drop the  term $ \beta v u_{xx } $  in (1.3),  the aforementioned violation of the law of gravity  will not occur.  It appears that  in the context of (1.3)-(1.4) the presence of this 2-nd order term $ \beta v u_{xx } $   principally changes the  ``mathematical nature'' of the problem  at hand, because it  becomes the senior order term in (1.3). Respectively, one cannot expect  its contribution to be  always  ``small''. 

In  Section 3 we show that  the 2-nd equation in the 1-$D$ model (1.3)-(1.4)  may  generate certain problems with the expected ``physical accuracy''  when the standing layer has corners.
To address this issue,  in {\em Section 4}, we  introduce a 1-$D$ {\em discrete} model for the dynamic granular matter formation described by  the system of {\em difference equations}, i.e.,  instead of the usual  partial  differential equations.  This approach seems to be novel to this field. Its  main idea is that it  will allow one to deal with the process at hand at the discretized  level and, hence,  to quantitatively distinguish between  the amounts of the rolling matter coming  down from the left and  from the right to  the {\em corner} points  of the standing layer. These  ``points of lack of differentiability'' are typical in the dynamic models when the standing layer is constantly re-created by the moving front of the rolling matter.  

Let us elaborate a little further  on the benefits of this  approach over the traditional pde approach in the  context of this  particular physical phenomenon. Clearly, the function representing the height of the standing layer becomes  non-differentiable at the corner points of the standing layer. The traditional pde approach (to modeling)  usually ``neglects'' the behavior of the process at hand at such ``isolated non-differentiable'' points (e.g.,  making use of the concept of generalized non-differentiable solutions). However, in the case of granular flow, the corners of standing layer (e.g., at the foundation of a sand pile) appear to be of major significance and cannot be overlooked  just as  ``unimportant isolated single points''. A possible way out here is the aforementioned difference equations approach.

We illustrate the advantages of our approach by  computational experiments in Section 7, which directly use the respective difference equations of Section 4 as a  computational algorithm. In particular, they show that the respective model (4.1)-(4.4) works  well   for any  slopes, regardless of  their  length, sign, steepness and thickness. 

In Section 5 we discuss the properties of model (4.1)-(4.4), which are used in  Section 6  to prove the  existence of  its   distributed limit solution  by  passing to the limit as the mesh size   tends to zero.
In  Section 8 we extend our 1-$D$ model to the case of two dimensions. Formal limits of our discrete models to those in the form of pde's, and their connection to existing models, are  also discussed.

{\bf 2. Diffusion term and the lack of compliance with  the  law of gravity  in model (1.3)-(1.4).}  In this section we will  re-visit  the phenomenon of ``geysers'' in model (1.3)-(1.4), discovered in \cite{Col} solely by numerical means. Our goal here is to investigate it by means of   formal analysis, which, in turn,  will provide us with an idea what can be done to modify model equations (1.3)-(1.4) to avoid this phenomenon.

Suppose that, initially, the granular matter was located strictly inside   of   (0, 1)  and that $  T > 0 $ is small enough to ensure that during the flow process  it would  remain within $ (0, 1)$. Thus, instead of (1.3)-(1.4) we will further deal with its homogeneous version in a bounded space domain:
$$
v_t = \beta (v u_x)_x - \gamma (\alpha - \mid u_x \mid)v, 
\eqno(2.1)$$
$$
u_t = \gamma (\alpha - \mid u_x \mid)v, \;\;\;\; (x,t) \in  Q_T = (0,1) \times (0, T),
\eqno(2.2)$$
$$
u(x,0) = u_o (x) \geq 0, \;\;  v(x,0) = v_o (x) \geq 0, \;\; \mid u_{ox} (x) \mid \; \leq \alpha, \;\;  \;\; x \in [0, 1],
$$
$$
u(x,t) \mid_{x = 0, 1} = 0, \;\;\;\; v(x,t) \mid_{x = 0, 1} = 0, \;\; t \in (0, T). 
$$
We formally assume below in our analysis  that the functions $ u, u_t, u_x, u_{xx}, v, v_t, v_x $  exist and are continuous in $ \bar{Q}_T$ and that there is no external source in model (1.3)-(1.4), i.e., $ f = 0$.

{\bf Lack of compliance with the law of gravity  in model (1.3)-(1.4)/(2.1)-(2.2).}
Suppose that, initially, the maximum value for $  \; (u_o (x) + v_o (x)) \;   $   is reached at point $ x_o, 0 < x_o <1, u_o (x_o) +  v_o (x_o) >  0$. Hence, $  v_{ox} + u_{ox} = 0$, or $  v_{ox}  = - u_{ox}$, and $  (v_o +  u_o)_{xx} \leq 0 \;\; {\rm at} \;\;  x_o$. In turn,  (2.1) yields that
$$
(u+v)_t = -\beta u^2_{ox} + \beta v_o u_{oxx}  \;\; {\rm at} \; (x_o, 0).
\eqno(2.3)$$
Now note that   the expression on the right in (2.3)  {\em can} be positive   when  the term  $ u_{oxx} $ is positive at $ x = x_o$. In this case, the  time-derivative of the height of granular matter   $ (u+v)_t \mid_{(x_o, t)}  $ will be positive at and  immediately after $ t = 0$, and thus it  will {\em have to increase} after $ t=0$ to exceed its  maximum initial value  $  \; (u_o (x_o) + v_o (x_o))$. In other words, model (2.1)-(2.2) allows the rolling matter to roll upward, ``defying''  the gravity law. 

{\bf Example 2.1: Geyser in  a cavity.} This  example is a version of the situation investigated  in \cite{Col} numerically. Suppose $ u_{o} (x) $ has a cavity near the maximum point $ x_0$ of $ u_o (x) + v_o(x)$ as shown in Fig. 1.

\setlength{\unitlength}{1mm}
\begin{picture}(150,55)(0,0)
\linethickness{1pt}
\put (2,15){\vector(1,0){85}} 

\put (10,10){\vector(0,1){35}}

\qbezier(10,15)(20,23)(28,21)
\qbezier(28,21)(34,20)(45,21)
\qbezier(45,21)(55,23)(63,15)

\put(35,10){\makebox(0,0)[b]{$x_{0}$}}

\put(35,22.3){\makebox(0,0)[b]{$v_{0}$}}

\put(35,15){\line(0,1){2}}

\put (24,21.5){\line(1,0){25}}

\put(8,40){\makebox(0,0)[b]{$u$}}
\put(85,10){\makebox(0,0)[b]{$x$}}

\put(8,10){\makebox(0,0)[b]{$0$}}
\put(63,15){\line(0,1){2}}
\put(63,10){\makebox(0,0)[b]{$1$}}
\put(40,15.7){\makebox(0,0)[b]{$u_o (x)$}}

\put(21.5,15){\line(0,1){2}}
\put(21,10){\makebox(0,0)[b]{$0.25$}}
\put(53,15){\line(0,1){2}}
\put(53,10){\makebox(0,0)[b]{$0.75$}}
\put(10,22){\line(1,0){2}}
\put(7,19){\makebox(0,0)[b]{$C_o$}}

\put(50,0){\makebox(0,10)[b]{${\rm Figure \; 1: \; ``Geyser\; situation''  \; in \; a  \; cavity.}$}}

\end{picture}

It is assumed that the maximum value $ C_o$ of the  function $ u_o (x) + v_o(x), x \in [0, 1]$  is reached {\em everywhere} between $ x = 0.25 $ and $0.75$ with 
the rolling layer  to be  filling the cavity  (e.g. sand ``pushed/dropped'' into a cavity). Then: 
$$
u_o (x) + v_o(x) = C_o, \;\;\;\; v_o (x)  = - u_o (x) + C_o, \;\; \;\; x \in [0.25, 0.75].
$$
We assume that $ u_o (x)$ is strictly convex at $ x_o$ and  $ u_o (x_o) $ is the point of local minimum for $ u_o (x)$, i.e., $ u_{oxx}(x_0) > 0$. Hence,  $ v_o (x)$ is strictly concave at $ x_o$, and $ v_o (x_o) $ is the point of local maximum for $ v_o (x)$. Respectively 
$$ 
v_{ox} =  u_{ox} = 0, \;\;\;\; (v_o +  u_o)_{xx} = 0 \;\; {\rm at} \;\;  x_o.
$$
Therefore, in (2.3): 
$$
(u+v)_t = \;   \beta v_o u_{oxx} > 0 \;\; {\rm at} \;  x= x_o, t = 0.
\eqno(2.3)^*$$
 Hence,   the value of $(v (x,t) + u (x,t)) $ will increase as $ t $ increases at $ x_o$ from its initial highest point, i.e., generating some sort of a ``geyser''  near this point.

{\bf Example 2.2: Geyser created by a ledge.} In this example we take only {\em  slightly more} than  the right half of the cavity to form a ``ledge'' on the slope of the standing layer as shown on Figure 2.  The same  argument applies to show  that one has a geyser at $ x_o$.

\setlength{\unitlength}{1mm}
\begin{picture}(150,55)(0,0)
\linethickness{1pt}
\put (2,15){\vector(1,0){85}} 

\put (10,10){\vector(0,1){35}}

\qbezier(10,15)(20,28.5)(28,31)
\qbezier(28,31)(36,30)(40,45)
\qbezier(40,45)(50,30)(63,15)

\qbezier(22,27.5)(22,31.5)(30,32.2)

\put(35,10){\makebox(0,0)[b]{$x_{0}$}}

\put(30,33.3){\makebox(0,0)[b]{$v_{0}$}}

\put(35,15){\line(0,1){2}}

\put (29,32){\line(1,0){4}}

\put(8,40){\makebox(0,0)[b]{$u$}}
\put(85,10){\makebox(0,0)[b]{$x$}}

\put(8,10){\makebox(0,0)[b]{$0$}}
\put(63,15){\line(0,1){2}}
\put(63,10){\makebox(0,0)[b]{$1$}}
\put(40,25.7){\makebox(0,0)[b]{$u_o (x)$}}

\put(20,15){\line(0,1){2}}
\put(20,10){\makebox(0,0)[b]{$0.25$}}
\put(53,15){\line(0,1){2}}
\put(53,10){\makebox(0,0)[b]{$0.75$}}
\put(10,32){\line(1,0){2}}
\put(7,29){\makebox(0,0)[b]{$C_o$}}

\put(50,0){\makebox(0,10)[b]{${\rm Figure \; 2: \; Geyser\; situation \; in \; a \; ``half- cavity''.}$}}

\end{picture}

{\bf Remark 2.1:  On ``thin'' and ``thick''   rolling slopes.}  We emphasize here that, due to (2.3),  the above situation does  not depend on the thickness of the rolling layer  relative to the standing layer. It rather depends on the shape (curvature) of the {\em surface of the standing layer}. In particular, in the above  counterexamples the depth of the cavity can be as   ``thin'' as one wishes on the surface of the standing layer of an arbitrary height, which implies that the respective rolling layer can  relatively be as ``shallow '' as one wishes.

{\bf Modification of equation (1.3)/(2.1).} We established in the above that  the diffusion term $ \beta v u_{xx} $ in model (1.3)-(1.4)/(2.1)-(2.2) can overpower the underlaying law of gravity. This questions the accuracy of this model, at least in some situations. If this term is {\em  dropped}, we will have the following model:
$$
v_t = \beta v_x u_x - \gamma (\alpha - \mid u_x \mid)v,
\eqno(2.4)$$
$$
u_t = \gamma (\alpha - \mid u_x \mid)v \;\;\;\; {\rm in} \;\;\;\; {\rm in} \; Q_T = (0,1) \times (0, T),
\eqno(2.5)$$
$$
u(x,0) = u_o (x) \geq 0, \;\;  v(x,0) = v_o (x) \geq 0, \;\; \mid u_{ox} \mid \; \leq \alpha \;\;  \;\; x \in [0, 1],
$$
$$
u(x,t) \mid_{x = 0, 1} = 0, \;\;\;\; v(x,t) \mid_{x = 0, 1} = 0 \;\; t \in (0, T). 
$$

Let us show that   model (2.4)-(2.5)  fully complies with the law of gravity.
The physical interpretation of  equation (2.4) is straightforward.
Let, say, $ u_x > 0, \;  u_x v_x > 0 $ at some $(x_0,t_0)$ as shown here:

\setlength{\unitlength}{1mm}
\begin{picture}(100,50)(0,0)
\linethickness{1pt}

\put (62,31){\vector(0,1){7}}
\put(62,39){\makebox(0,0)[b]{$(u + v)$}}

\put (82,36){\vector(-2,-1){7}}

\put(20,10){\line(2,1){70}}
\put(20,10){\line(3,1){70}}
\put(10,10){\line(1,0){70}}

\put(62,10){\line(0,1){2}}

\put(65,20){\makebox(0,0)[b]{$u$}}

\put(70,32){\makebox(0,0)[b]{$v$}}

\put(65,12){\makebox(0,0)[b]{$x_0$}}

\put(50,0){\makebox(0,10)[b]{${\rm Figure \; 3: u_x, u_x v_x > 0 \; {\rm at} \;   (x_0, t_0).}$}}

\end{picture}

\par\noindent
On Fig. 3 the standing layer is rising at $ x_0$ at time $t_0$, and {\em relative to it},  the rolling layer is rising as well. Hence, the granular matter in the latter  rolls down (to the left) and this {\em will} increase the value of $ \; (u + v) \; $ at this point, as confirmed by (2.4),  rewritten as follows:
$$
(v + u) \mid_{(x_0,t)} \;  = \; (v + u) \mid_{(x_0,t_0)} \; + \; \int_{t_0}^t v_x (x_0, \tau) u_x (x_0, \tau) d \tau.
$$
The next picture illustrates  the case when $ u_x > 0, \;  u_x v_x < 0 $ at some $(x_0,t_0)$.

\setlength{\unitlength}{1mm}
\begin{picture}(100,50)(0,0)
\linethickness{1pt}

\put (60,34){\vector(0,-1){7}}
\put(65,35){\makebox(0,0)[b]{$(u + v)$}}
\put (52,26){\vector(-2,-1){7}}

\put(80,30){\line(-5,1){30}}
\put(20,10){\line(3,1){70}}
\put(10,10){\line(1,0){70}}

\put(60,10){\line(0,1){2}}

\put(65,20){\makebox(0,0)[b]{$u$}}

\put(45,32){\makebox(0,0)[b]{$v$}}

\put(60,10){\line(0,1){2}}
\put(65,12){\makebox(0,0)[b]{$x_0$}}

\put(50,0){\makebox(0,10)[b]{${\rm Figure \; 4: u_x > 0, u_x v_x < 0 \; {\rm at} \;   (x_0, t_0).}$}}

\end{picture}

 {\bf  Proposition  2.1: Compliance of model (2.4)-(2.5) with the law of gravity.} }
 {\em The maximum value of  $ \; u (x,t) + v (x,t) \; $ for the model (2.4)-(2.5)  is attained  at $ t = t_0$.}

{\em Proof.}  The proof  is quite simple and makes use of the idea of the classical proof of the maximum principle for the linear 1-$d$ heat equation. 

We argue by contradiction and consider only the nontrivial case: $ v + u \not\equiv 0$ in $Q_T$. Let   the function $ (u(x,t)  + v (x,t) ) $ attain its  global maximum in $ [0,1] \times [0, T]$,  which is {\em strictly greater} than  $ \max_{x \in [0,1]} (v_o (x) + u_o (x)) $,  at some point $ x_o \in (0, 1)$ ({recall here  that $ (v  + u) \mid_{x =0, 1} = 0$) at a positive moment time. Then, by continuity,  the same is true for the function $ (v^*(x,t)  + u^*(x,t)) $, where  $ v^* = e^{-kt}v, \;\; u^* = e^{-kt}u $
and   $k>0$ is a ``small'' positive parameter. 
Now, let  the  just-mentioned ({\em non-zero}) global maximum for the function $ (v^*(x,t)  + u^* (x,t))$ be attained at  $(x^*,  t^*)$, where  $ x^* \in (0,1), t^* >0$.
If so, we should have
$$
e^{-kt} (v  + u)_x = (v^*  + u^*)_x   =   0, \;\; (v^*  + u^*)_t  \geq   0 \;\; {\rm at} \; (x^*,  t^*) 
$$
({we  use  the  left-hand side  time-derivative  if $\; t^* = T$). However, due to (2.4):
$$
(u^*+v^*)_t = -ke^{-kt}(v + u) - \beta e^{-kt} v^2_{x} < 0 \;\; {\rm at} \; (x^*,  t^*).
$$
Hence,  we arrived at contradiction.
End of proof.

 {\bf Remark 2.2: Comparison of Proposition 2.1 with  the classical maximum principle for parabolic pde's.} 
There is a resemblance  between Proposition 2.1 and the classical  maximum principle for the {\em linear} 1-$d$ heat equation, which states that the extreme values of the temperature  within a spatial domain are attained either at the initial moment or on the boundary of this domain. This principle reflects the main underlying physical law of the heat transfer. However, the heat equation is a {\em scalar linear} pde,  while Propositiuon 2.1 deals with a 2-$D$ {\em nonlinear system} of pde's. Thus, it makes it  a novelty how one can setup a ``maximum principle'' for system like (2.4)-(2.5). Note that Proposition 2.1 does {\em not} deal with the separate maximum values of  functions $ u $ and $ v$. 
}

{\bf 3.  Behavior of model (1.3)-(1.4)/(2.1)-(2.2) near the corners in the standing layer created by the  critical slope(s).} Our goal here is to explain the reasons for our further modifications of this model as suggested below.

Equation  (2.2)/(1.4) does not  allow  the standing layer to  grow at the points where the slope is critical, i.e., when $ \mid u_x \mid = \alpha$ and, thus,  $ u_t = 0$. 
Consider, e.g.,  the  case when  the base of 
the standing layer is extending due to influx of the rolling matter, as illustrated by  Figures 5 and 6:

\setlength{\unitlength}{1mm}
\begin{picture}(100,50)(0,0)
\linethickness{1pt}

\put (82,38){\vector(-2,-1){7}}

\put(20,10){\line(2,1){70}}
\put(20,10){\line(5,2){70}}
\put(0,10){\line(1,0){90}}

\put(68,18){\makebox(0,0)[b]{$u (x,t) = \; {\rm constant}$}}

\put(60,34){\makebox(0,0)[b]{$v (x,t) $}}

\put(21,10){\line(0,1){2}}
\put(19,12){\makebox(0,0)[b]{$x_*$}}

\put(50,0){\makebox(0,10)[b]{${\rm Figure \; 5: \; The \; case \; of \; critical \; slope.}$}}

\end{picture}

\par\noindent
Let assume that the slope of the standing layer  $u$ on Fig. 5  is critical ($ u_x = \alpha$)  on the right of the corner point $x_*$, i.e.,  where it meets  the horizontal base. Then, due to (1.4)/(2.2),  
$$
u_t (x,t) = 0, \;\;\; u(x,t) = u(x,0) = u_o (x),  \;\;\; t > 0, \;\; x > x_*,
\eqno(3.1)$$
and, thus,  all the material in the rolling layer will have to roll to the left of this corner point,  leaving the  standing layer unchanged (see also Remark 3.1below).
However, the physical expectations here are associated with  a build-up of the standing layer {\em on both sides of $x_*$}, e.g., as shown on Fig.~6:

\setlength{\unitlength}{1mm}
\begin{picture}(100,50)(0,0)
\linethickness{1pt}

\put(10,10){\line(5,1){38}}

\put(48,18){\line(2,1){28}}
\put(0,10){\line(1,0){90}}

\put(50,11){\makebox(0,0)[b]{$u(x,t) $}}

\put(34,10){\line(0,1){2}}
\put(36,12){\makebox(0,0)[b]{$x_*$}}

\put(55,0){\makebox(10,10)[b]{${\rm Figure \; 6: \; Expected \; evolution \; of \; the \; standing \; layer \; from  \; Fig.  \, 5.}$}}

\end{picture}

Fig. 6 is not compatible with  (1.4)/(2.2).  Thus, this equation is  not accurate  in this case, i.e.,   when there is a transition from the critical slope to a lesser one.

 {\bf Remark 3.1.} In Figs. 5 the standing layer has a {\em nondiffrentiable} slope at $ x_*$. Hence, model (1.3)-(1.4)/(2.1)-(2.2) cannot be applied  in the classical sense. If, to deal with it (say, numerically), we assume that
$$ 
u_x (x_*, 0) = \frac{1}{2} (0 + \alpha),
$$
i.e., take the average of slopes from the right and from the left, this will only result in the growth of $ u (x,t)$ at $ x_*$ due to equation (1.4)/(2.2) (this growth at a single point does not make ``physical'' sense). However, on the right of $ x_*$ the slope will remain critical and hence still cannot grow into a physically  expected formation as on  Figure 6.

{\bf Cavity.} Similar to the above, in the case when the standing layer has a cavity as on Figure 7(b) below with both slopes converging to the same point to be critical, equations (1.4)/(2.2) will not allow these slopes to grow when there is a rolling matter on either of them. 
We propose a  way out in this situation in the next section (see also numerical examples 7.1-2 in Section 7). 

{\bf ``Abrupt  halt and full conversion'' of the rolling matter.} If there is no inertia, then, when the rolling matter arrives to the {\em lowest}  corner point of the cavity, it  should immediately stop and instantaneously be {\em fully} converted into the standing layer at this point. This phenomenon is impossible to describe by means of pde's like in (1.3)-(1.4) (one may consider, e.g., $\delta$-functions for that).
 
{\bf 4. 1-$D$ difference equations model.}  The discussion in Section 3 indicates that the situations when the standing layer has corners create serious problems for modeling by means of pde's. Mathematically, these situations  are associated with the spatial points  where the 1-st spatial derivative $ u_x$ does not exist. In our opinion, this justifies attempts to try to  look for different types  of  model equations  which do not deal with differentiable (both in the classical and generalized sense) functions. 

\underline{Our approach} below,  to address the just-outlined  issue, deals with an  attempt to {\em convert  the above-discussed pde model}  (2.4)-(2.5) into an associated  {\em discrete    model}, which would intrinsically allow us to use one-sided derivatives  to distinguish between  the amounts of the rolling matter coming  down from the left and  from the right to  the corner points  of the standing layer (or rolling down away from them).  In the process of  construction of  this model, we intend to focus on its  physical interpretation and to discuss all possible cases of mutual orientations of slopes at hand.

Let us split the interval [0, 1] into $n$ equal segments of size $\Delta_{n}  x$, and the time-interval $[0, T]$ into $m$ equal segments of size 
$ \Delta_m t$.

{\bf Remark 4.1.} In what follows,  the size of time-step will be allowed to  change to  smaller  values before we reached the moment $T$. Nonetheless, for simplicity of notations we will use the same $m$ below. The space-step remains constant for the given system. (One, of course, can select a different mesh strategy.)

The solution of the approximate system is represented by the collection of values denoted by $ \{ u^{(n, m)}_{i, j}, v^{(n, m)}_{i, j}, i = 0, \ldots, n, j = 0, \ldots \}$. These values define the   {\em piecewise linear}  approximate solution, {\em  denoted for the standing layer  by $  u^{(n, m)} (x, t)$ and for the rolling layer by $  v^{(n, m)} (x, t), x \in [0,1], t \in [0, T]$}, to the granular matter formation process at hand for the respective nodes $ \{(x_i, t_j), i = 1, \ldots, n, j = 0, \ldots\}$. For example,
$$
u^{(n, m)}_{0, 0}  = u^{(n, m)} (0, 0), \;\; u^{(n, m)}_{1, 1}  = u^{(n, m)} (\Delta_{n}  x,  \Delta_{m}  t), \;\; {\rm and \; so \; forth.}
$$

The physical idea behind the propagation of rolling matter, exploited  in this paper,  is that it can move {\em only downward}  due to the force of gravity. Therefore, the changes happening to the combined rolling and standing layers  at point $x_i = i \Delta_n  x$ are to be the result only of the {\em matter (a) arriving down to it from the adjacent rising slope(s)}  and/or {\em (b) leaving down along the adjacent falling slope(s)}.

{\bf Remark 4.2.} The above means that, in what follows, we do not take into account other possible motions of the granular matter such as, e.g., of   stochastic nature or due to inertia, etc.

Note that, unlike  models (2.1)-(2.2) and (2.4)-(2.5),  in the case of difference equations approach we will  deal with {\em the left- and right-hands sides derivatives} of our piecewise linear approximate functions  $  u^{(n, m)} (x, t)$. 
Three principal cases are possible near the point $x_i$ ({\em case (7(a)  has also its symmetric double about the vertical line passing through $ x_i$}):

\setlength{\unitlength}{1mm}
\begin{picture}(100,40)(0,0)
\linethickness{1pt}

\put(42,10){\line(0,1){2}}
\put(21,10){\line(0,1){2}}
\put(60,10){\line(0,1){2}}

\put(42,17.2){\line(2,1){18}}
\put(20,10){\line(3,1){22}}
\put(10,10){\line(1,0){70}}

\put(40,12){\makebox(0,0)[b]{$x_{i}$}}
\put(65,12){\makebox(0,0)[b]{$x_{i+1}$}}
\put(14,12){\makebox(0,0)[b]{$x_{i-1}$}}

\put(60,0){\makebox(0,10)[b]{${\rm Figure \; 7(a): \; Slope \; is \;   rising \; to \; the \; right \; of \; x_i \; only. \; Both \; slopes  \; can \; be \; zero.}$}}

\end{picture}

\setlength{\unitlength}{1mm}
\begin{picture}(100,40)(0,0)
\linethickness{1pt}

\put(43,10){\line(0,1){2}}
\put(21,10){\line(0,1){2}}
\put(60,10){\line(0,1){2}}

\put(43.5,16.8){\line(2,1){16.5}}
\put(21.5,24){\line(3,-1){21.9}}

\put(10,10){\line(1,0){70}}

\put(44,12){\makebox(0,0)[b]{$x_{i}$}}
\put(65,12){\makebox(0,0)[b]{$x_{i+1}$}}
\put(15,12){\makebox(0,0)[b]{$x_{i-1}$}}

\put(55,0){\makebox(0,10)[b]{${\rm Figure \; 7(b):  \; A \; ``cavity'' \; near \; x_i\; with \; slopes \;  rising \; on \; both \; sides \; of \; x_i.}$}}

\end{picture}

\setlength{\unitlength}{1mm}
\begin{picture}(100,40)(0,0)
\linethickness{1pt}

\put(43,10){\line(0,1){2}}
\put(21,10){\line(0,1){2}}
\put(60,10){\line(0,1){2}}

\put(43.5, 32.5){\line(2,-1){16.5}}
\put(21.5,18){\line(3,2){21.9}}

\put(10,10){\line(1,0){70}}

\put(44,12){\makebox(0,0)[b]{$x_{i}$}}
\put(65,12){\makebox(0,0)[b]{$x_{i+1}$}}
\put(15,12){\makebox(0,0)[b]{$x_{i-1}$}}

\put(55,0){\makebox(0,10)[b]{${\rm Figure \; 7 ({\rm c}):  \; A ``vertex/peak'' \; at\; x_i\; with \; slopes \;  falling \; on \; both \; sides \; of \; x_i.}$}}

\end{picture}

\bigskip

To approximate the time derivative at point $ (x_i, t_j)$, we use the standard forward approximation:
$$
u_t \approx \frac{u^{(n, m)} (x_i, t_{j+1}) - u^{(n, m)} (x_i, t_{j})}{\Delta_m t}.
$$
To approximate the spatial  derivative at point $ (x_i, t_j) $, we use respectively
 the  derivatives of   $  u^{(n, m)} (x, t)$: 
$$
u_{x+}^{(n, m)} (x_i, t_j)  \; = \; \frac{u^{(n, m)} (x_{i+1}, t_{j}) - u^{(n, m)} (x_i, t_{j})}{+\Delta_n x},   
$$
$$
u_{x-}^{(n, m)} (x_i, t_j)  \; = \; \frac{u^{(n, m)} (x_{i-1}, t_{j}) - u^{(n, m)} (x_i, t_{j})}{-\Delta_n x}. 
$$
Our equations (4.1)-(4.2) below are derived based on  the discussions around Figures 3-7(a-c).

{\bf Difference equations for the  standing layer.} In this case we assume that the increase of the height of the standing layer is defined  by  the slope(s) of the standing layer  {\em below} the given grid point as the one(s)  determining:
 \begin{itemize}
  \item
  how much of  the rolling matter {\em available  at this point will \underline{stay} there} 
  \item
  and {\em how much of it will \underline{roll down}}, namely:  
  \end{itemize} 
  $$
u^{(n, m)} (x_i, t_{j+1})   = u^{(n, m)} (x_i, t_{j}) 
$$
$$
+\;  \;   
\gamma \Delta_m t  \; \left\{ \begin{array}{ll}
  \;  (\alpha  - \mid u_{x-}^{(n, m)} (x_i, t_j)\mid)  v^{(n, m)} (x_i, t_{j})   &  {\rm for} \; u_{x+}^{(n, m)} (x_i, t_j)  \geq 0, \;  \; u_{x-}^{(n, m)} (x_i, t_j)  > 0\;  \\ 
 \;  (\alpha  - \mid u_{x+}^{(n, m)} (x_i, t_j)\mid)  v^{(n, m)} (x_i, t_{j})   &  {\rm for \;}  u_{x-}^{(n, m)} (x_i, t_j)  \leq 0, \;  \; u_{x+}^{(n, m)} (x_i, t_j)  < 0  \\ 
\frac{1}{\gamma \Delta_m t}   v^{(n, m)} (x_i, t_{j})  &  {\rm for \;  }   
u_{x+}^{(n, m)} (x_i, t_j)  \geq 0, \;  \; u_{x-}^{(n, m)} (x_i, t_j)  = 0, \\
\frac{1}{\gamma \Delta_m t}   v^{(n, m)} (x_i, t_{j})  &  {\rm for }    \; u_{x-}^{(n, m)} (x_i, t_j)  \leq 0, \;  \; u_{x+}^{(n, m)} (x_i, t_j)  = 0, \\
\frac{1}{\gamma \Delta_m t}    v^{(n, m)} (x_i, t_{j})  &  {\rm for \;}     \; u_{x+}^{(n, m)} (x_i, t_j)  \geq 0, \;  \; u_{x-}^{(n, m)} (x_i, t_j)  \leq 0, \\
 r_{i,j,_-}   (\alpha 
  - \mid u_{x-}^{(n, m)} (x_i, t_j)\mid)  v^{(n, m)} (x_i, t_{j})  \\
+  \;  r_{i,j,_+} (\alpha  
- \mid u_{x+}^{(n, m)} (x_i, t_j)\mid)  v^{(n, m)} (x_i, t_{j})  &  {\rm for \; }    \; u_{x+}^{(n, m)} (x_i, t_j)  < 0,  \;  \; u_{x-}^{(n, m)} (x_i, t_j) >  0 .
  \end{array}
\right.
\eqno(4.1)$$
Note that formulas (4.1)  will not immediately  increase the slopes  beyond the critical value both to the left and to the right of point  $x_i$ for Fig. 7(a) and/or its symmetric counterpart -- if  $\Delta_mt$  is sufficiently small (see Remark 4.1). When the lower slope is zero or like on Fig. 7(b), there is no lower slope for rolling down from $ x_i$, and therefore all  the rolling matter at this point becomes an addition to the standing layer at this point. (Let us remind the reader that in our model we do not take inertia into consideration.) 
In the case of  Fig. 7(${\rm c}$) we assume that the rolling matter will all roll down along the steepest slope as defined by the following choice of splitting coefficients $r_{i,j,_\pm} $: 
$$
r_{i,j,_-} \; = \; 
\left\{ \begin{array}{ll}
1 &  {\rm if } \; \mid u_{x-}^{(n, m)} (x_i, t_j)\mid  >  \mid u_{x+}^{(n, m)}  (x_i, t_j)\mid\;  \\ 
0 &  {\rm if } \; \mid u_{x-}^{(n, m)} (x_i, t_j)\mid  <  \mid u_{x+}^{(n, m)}  (x_i, t_j)\mid\;  \\ 
\frac{1}{2} &  {\rm if } \; \mid u_{x-}^{(n, m)} (x_i, t_j)\mid  =  \mid u_{x+}^{(n, m)}  (x_i, t_j)\mid,
  \end{array}
\right.
$$
while  $ r_{i,j,_+}$  is defined symmetrically. This strategy   will ensure, in particular,  that either slope, on the left and on the right of $ x_i$, will not exceed the critical value. Namely, if at least one of the aforementioned slopes is critical, then the standing layer will not increase at $ x_i$. 
 
Denote 
$$
v_*^{(n, m)} (x_i, t_{j})=  v^{(n, m)} (x_i, t_{j}) -\left(u^{(n, m)} (x_i, t_{j+1})  - u^{(n, m)} (x_i, t_{j}) \right), 
$$
where (and in (4.1)) the term in the large parenthesis describes the increase of the height of the standing layer during the time-interval $ (t_j, t_{j+1})$ due to the contribution from the rolling  layer available at time $ t_j$. Thus, $
v_*^{(n, m)} (x_i, t_{j}) $ describes the part of  the rolling layer at $x_i$ which will leave  this point rolling down a respective slope (unless we have the situation as on Fig. 7(${\rm b}$)) after $ t = t_j$.  In other words, we assume that the contribution of the rolling layer to the standing layer on the interval $ [t_j, t_{j+1}]$ occurs at time $t_j$. Respectively, in equations (4.2) for the rolling layer,  we will use the following notations:
$$
v_{*x+}^{(n, m)} (x_i, t_j)  \; = \; \frac{v_*^{(n, m)} (x_{i+1}, t_{j}) - v_*^{(n, m)} (x_i, t_{j})}{+\Delta_n x},   \;\;\;\; 
$$
$$v_{*x-}^{(n, m)} (x_i, t_j)  \; = \; \frac{v_*^{(n, m)} (x_{i-1}, t_{j}) - v_*^{(n, m)} (x_i, t_{j})}{-\Delta_n x}.
$$

 {\bf Difference equations for the  rolling layer}: 
$$
v^{(n, m)} (x_i, t_{j+1})=  v_*^{(n, m)} (x_i, t_{j}) \; +\; \Delta_m t \;   \beta \; {\mathbf F},
\eqno(4.2)$$
where (see also Fig. 8 and explanations to it):
  $$
{\mathbf F}  =   r_{i+1,j+1,_-}  u_{x+}^{(n, m)} (x_i, t_{j + 1}) v_{*x+}^{(n, m)} (x_i, t_j)
 $$
 for the case like on Fig.  7(a)  when $ \;   u_{x+}^{(n, m)} (x_i, t_{j + 1})  \geq 0,   \; u_{x-}^{(n, m)} (x_i, t_{j + 1}) \geq 0$;
$$
{\mathbf F} = 
 r_{i-1,j+1,_+}  u_{x-}^{(n, m)} (x_i, t_{j + 1}) v_{*x-}^{(n, m)} (x_i, t_j)
 $$
for $  u_{x-}^{(n, m)} (x_i, t_{j + 1})  \leq 0,    \; u_{x+}^{(n, m)} (x_i, t_{j + 1})  \leq  0$;
$$
 {\mathbf F}=  r_{i+1,j+1,_-} u_{x+}^{(n, m)} (x_i, t_{j + 1}) v_{*x+}^{(n, m)} (x_i, t_{j + 1})  
+  r_{i-1,j+1,_+} u_{x-}^{(n, m)} (x_i, t_{j + 1}) v_{*x-}^{(n, m)} (x_i, t_{j + 1}) 
$$
for the case like on Fig.  7(b);
$$
 {\mathbf F} =  - \frac{1}{\Delta_m t \;   \beta}v_*^{(n, m)} (x_i, t_{j})
 $$
 for Fig.  7(${\rm c}$)  when     $  u_{x+}^{(n, m)} (x_i, t_{j + 1})  \leq 0, \;  \; u_{x-}^{(n, m)} (x_i, t_{j + 1})  \geq 0 $,  i.e., when there is  no influx of  rolling  matter  (all  the  prior  matter  will  leave  the  peak).

\setlength{\unitlength}{1mm}
\begin{picture}(100,50)(0,0)
\linethickness{1pt}

\put(43,10){\line(0,1){2}}
\put(21,10){\line(0,1){2}}
\put(60,10){\line(0,1){2}}

\put(5,40.5){\line(1,-1){16.5}}

\put(43.5,16.8){\line(2,1){16.5}}
\put(21.5,24){\line(3,-1){21.9}}

\put(60,25){\line(3,-1){21.9}}

\put(10,10){\line(1,0){70}}

\put(44,12){\makebox(0,0)[b]{$x_{i}$}}
\put(65,12){\makebox(0,0)[b]{$x_{i+1}$}}
\put(15,12){\makebox(0,0)[b]{$x_{i-1}$}}

\put(68,3){\makebox(0,10)[b]{${\rm Figure \; 8:  \; A \; ``cavity'' \; near \; x_i\; with \; slopes \;  rising \; on \; both \; sides \; of \; x_i.}$}}
\put(68,-2){\makebox(0,10)[b]{${\rm At  \; x_{i+1}  \; we \; have \; an  \; ``above \; peak'', \; no ``peak'' \; at \; x_{i-1}}$.}}

\end{picture}

\bigskip
\par\noindent
Figure 8 illustrates the idea of the choice of splitting coefficients in (4.2). To describe the forthcoming changes  in the height of the rolling layer  at  point $x_i$, we rely only on the slopes  {\em above}   the respective grid point,  i.e.,  where the rolling matter is coming from.

{\bf The initial and boundary conditions} for equations (4.1)-(4.2)  are defined similar to those in  (2.1)-(2.2):
$$
 u^{(n, m)} (0, t_j)  = u^{(n, m)} (1, t_j) = 0, \;\;\;\; j = 0, 1, \ldots,
$$
$$
u^{(n, m)} (x_i, 0) = u_0 (x_i) \; \geq 0, \;\; v^{(n, m)} (x_i, 0) = v_0 (x_i, 0) \; \geq 0, \;\; i = 0, \ldots, n,
\eqno(4.3)$$
$$
\mid u_{x\pm}^{(n, m)} (x_i, 0)   \mid  \; \leq \;\alpha, \;\; i = 1, \ldots, n-1. 
\eqno(4.4)$$

 {\bf 5. Properties of solutions to the difference equations.} Let us  remind the reader that we assumed in the above that we can  regulate the size of  time step in (4.1)-(4.2) to preserve  conditions (4.3) and (4.4) (see Remark 4.1).

{\bf Property 5.1: Compliance of model (4.1)-(4.4) with the law of gravity.} Similar to Proposition 2.1, 
the following inequality holds {\em with proper selection of the time-steps} in (4.1)-(4.2):
 $$
\max_{i = 1, \ldots, n; j = 1, \ldots } \{u^{(n, m)} (x_i, t_{j}) + v^{(n, m)} (x_i, t_{j}) \} \; 
$$
$$
\leq \; \max_{i = 1, \ldots, n} \{u^{(n, m)} (x_i, 0) + v^{(n, m)} (x_i, 0)\}
\eqno(5.1)$$
 Indeed, it follows from  (4.2) (see the last line)  that  in order to have
$$
u^{(n, m)} (x_{io}, \Delta_m t) + v^{(n, m)} (x_{io}, \Delta_mt) \; \geq \;  u^{(n, m)} (x_{io}, 0) + v^{(n, m)} (x_{io}, 0)
$$
for some  $ \Delta_m t$, one of the  slopes $u_{x\pm}^{(n, m)} (x_i, \Delta_m t) $ should rise away from the point $ x_{io}$ and the respective slope $ v_{*x \pm}^{(n, m)} (x_{io}, 0)$ should be of the same  sign too. This is impossible, since
$$
u^{(n, m)} (x_{io}, \Delta_m t) + v_*^{(n, m)} (x_{io}, 0) \; = \;  u^{(n, m)} (x_{io}, 0) + v^{(n, m)} (x_{io}, 0)
$$
and the right-hand side is the top vertex for the combined sum of layers at $ t = 0$.

Hence,  a possible inequality contradicting to (5.1) at time $t$ may arise only at some $ x_i \neq x_{i_0}$ with an initially strictly lower sum of two layers. Therefore, due to the fact that these are finitely many, we can select, if necessary from now on a new, smaller step size $ \Delta_{m^*} t$ which will guarantee (5.1) at time $ t = \Delta_{m}^*  t$ with the equality at  such spatial node point(s).

{\bf Property 5.2: The critical slope restriction.}  At time $ t = 0$ we have necessary conditions for that as given by (4.4). The explanations after (4.1) show that  this restriction cannot be violated immediately.
If, however,  $ u_{x+}^{(n, m)} (x_i, t_j)$ or $u_{x-}^{(n, m)} (x_i, t_j) $ will exceed the critical value at some moment $ t_j$, we are to return to the prior time-layer and to select a new smaller value for the time-step  which will not result in the violation of the critical slope requirement.  This process can go on till  $ u_{x}^{(n, m)} (x_i, t_j)$ had reached the critical slope everywhere and thus we have reached the saturated (static) solution.

{\bf Property 5.3: Nonnegativity of solutions.}  Equations (4.1) imply that in order for some of $u^{(n, m)} (x_{i }, t_{j+1})$'s become negative, some of $v^{(n, m)} (x_{i}, t_j)$ must become negative first. Let, e.g., at some point $(x_i, t_j)$ ($t_j $ can be zero) the value of  $v^{(n, m)} (x_{i}, t_j)$ becomes zero for the first time (for suitably adjusted  $ \Delta_m t$), while  $u^{(n, m)} (x_{i }, t_j) \geq 0$. In this case the sum of two layers at time $ t_{j + 1}$ will not decrease at $ x_{i}$.
Indeed, in this case the slopes $ v_{*x \pm}^{(n, m)} $ on either  side of $ x_i$ are either equal to zero or are rising {\em away} from it at time $t_j$. Hence,  the last term on the right in (4.2)  can be nonzero only if $ u_{x \pm}^{(n, m)} $ are also rising away from $ x_i$ at time $ t_{j + 1}$, because  (4.2) takes into account only the standing layer  slope(s) which lie above the respective grid point. Hence,  the last term on the right of (4.2) is non-negative in our case and therefore the sum of two layers cannot decrease at $x_i$.

Next, since $v_*^{(n, m)} (x_{i }, t_j) = v^{(n, m)} (x_{i }, t_j) = 0$, (4.1) means that 
$$
u^{(n, m)} (x_{i }, t_{j+1}) = u^{(n, m)} (x_{i }, t_j) \geq 0 \;\; {\rm at} \; x_i.
$$
Thus, for the aforementioned sum of two layers not to decrease at time $ t_j$,  the value of $v^{(n, m)} (x_{i}, t_{j+1})$ must remain nonnegative.

{\bf Property 5.4: Preservation of growth of $u^{(n, m)} (x_i, t_j)$ in time.} This follows immediately from (5.2) under Properties  5.2-3.

{\bf 6. Existence of a distributed solution to the dynamic granular  formation problem as a limit solution to model (4.1)-(4.4).} 
Note that the properties 5.1-4 mean that the sequence of functions $ \{ u^{(n, m)} (x, t)\}_{\Delta_n x, \Delta_m t} $ is uniformly bounded in  $ H_0^{1,0} (Q_T) =\{\phi \mid \phi, \phi_x \in L^2 (Q_T), \phi \mid_{x = 0, 1} = 0\}$, and is equicontinuous and uniformly bounded in $ C (\bar{Q}_T)$. Therefore, the Arzela-Ascoli Theorem yields that it contains a uniformly converging subsequence. Without loss of generality, we can say that
$$
u^{(n, m)} (x, t)  \; \rightarrow \;  u_* (x, t ) 
$$
as $ \;  \Delta_n x, \Delta_m t  \rightarrow 0+$   in $  \;  C (\bar{Q}_T)$  and weakly in $\; H^{1,0} (Q_T)$.
This limit function satisfies all the Properties 5.1, 5.3-4 in the continuous form. In turn, the slope restriction is satisfied in the ``generalized'' sense: 
the graph of function $ u_* (x, t) $ in $x $ for any $ t$ lies within the cone with center at point $ (x, u_* (x,t)) \in R^2$ and slopes $ \pm1$, while 
$ u_* (x,t) $ is non-decreasing in time for any $ x$.

Respectively,  the sequence of functions $ \{ v^{(n, m)} (x, t)\}_{\Delta_n x, \Delta_m t} $ is uniformly bounded in  $ L^\infty (Q_T)$, and thus  also  in $ L^2 (Q_T)$. Therefore, without loss of generality, we can say that
$$
v^{(n, m)} (x, t)  \; \rightarrow \;  v_* (x, t ) \;\;{\rm as} \; \Delta_n x, \Delta_m t  \rightarrow 0+ \;\; {\rm weakly \; in} \;  L^2 (Q_T).
$$

{\bf Connection to pde modeling.}  The difference equations (4.1)-(4.4) can be viewed as approximation of pde model  (2.4)-(2.5) at    points $ (x,t) $ where $ u $ and and $v$ are continuously differentiable and $ u_x \neq 0$.

{\bf 7. Model (4.1)-(4.4), computational strategy and examples.} In model (2.1)-(2.2) the second-order diffusion term is "symmetric", while the first-order convective term is "directional". For small values of $\beta$ the diffusion is ``empirically small'' compared to convection, i.e., in some situations the model exhibit strong propagation behavior for the most of the domain. In this case approximation of the convective terms by central difference leads to nonphysical oscillations \cite{Thomas}. On the other hand the proposed difference model (4.1)-(4.4) takes care of the aforementioned phenomena. The finite differences in the model are adaptively directional, taking into account the direction of propagation at each space step, thus allowing simulate accurately various configurations of standing and rolling matter. 
In  the two examples 7.1 and 7.2  illustrated respectively by Figures 9-11 and 12-14 we employ directly the difference equations (4.1)-(4.4) as a numerical algorithm, assuming uniform finite difference mesh in space and time and
$\alpha=1$, $\gamma=1$.

\setcounter{figure}{8}

\begin{center}
\begin{figure}
\includegraphics[width=7.5cm, height=4cm]{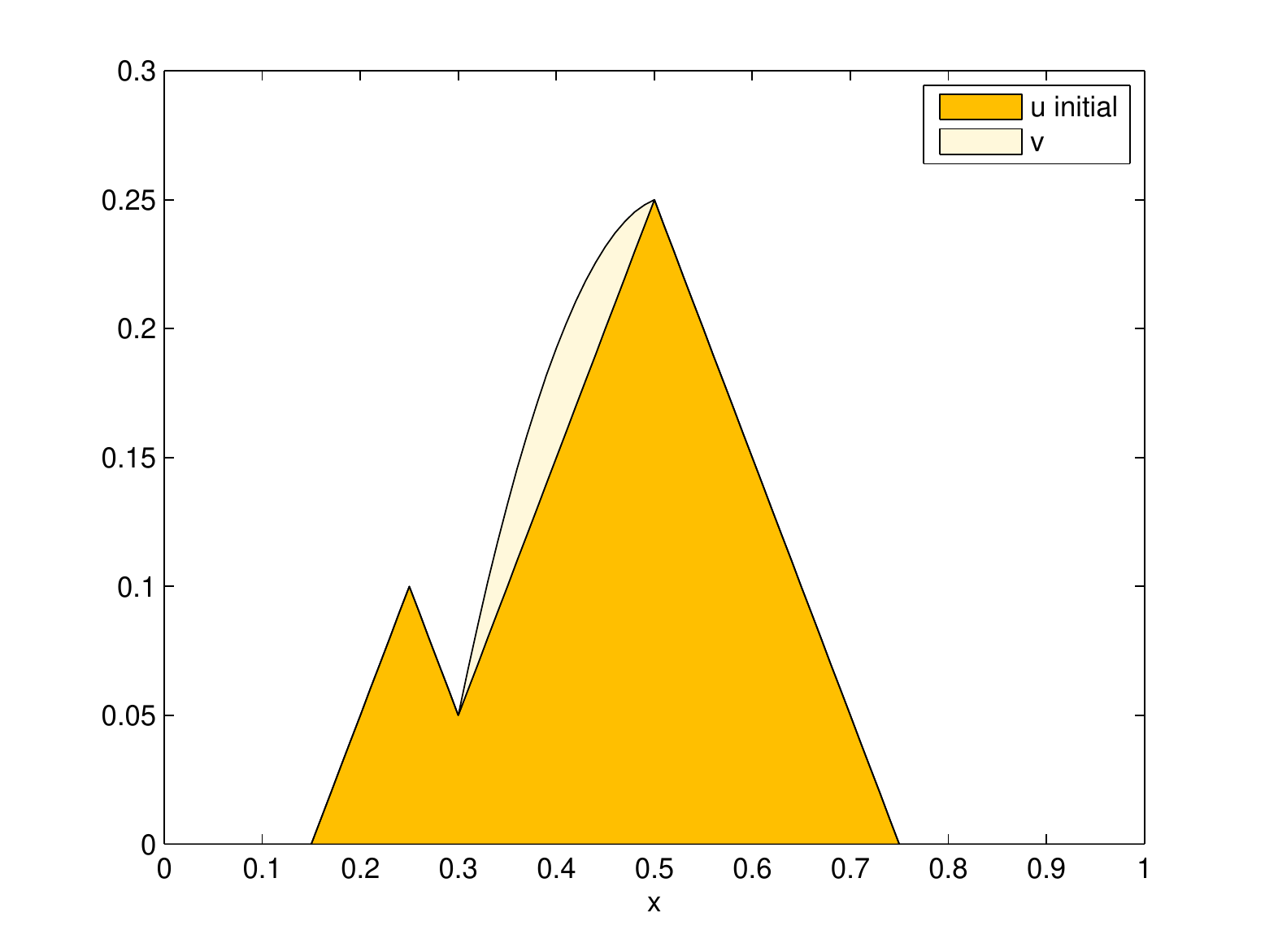}
\includegraphics[width=7.5cm, height=4cm]{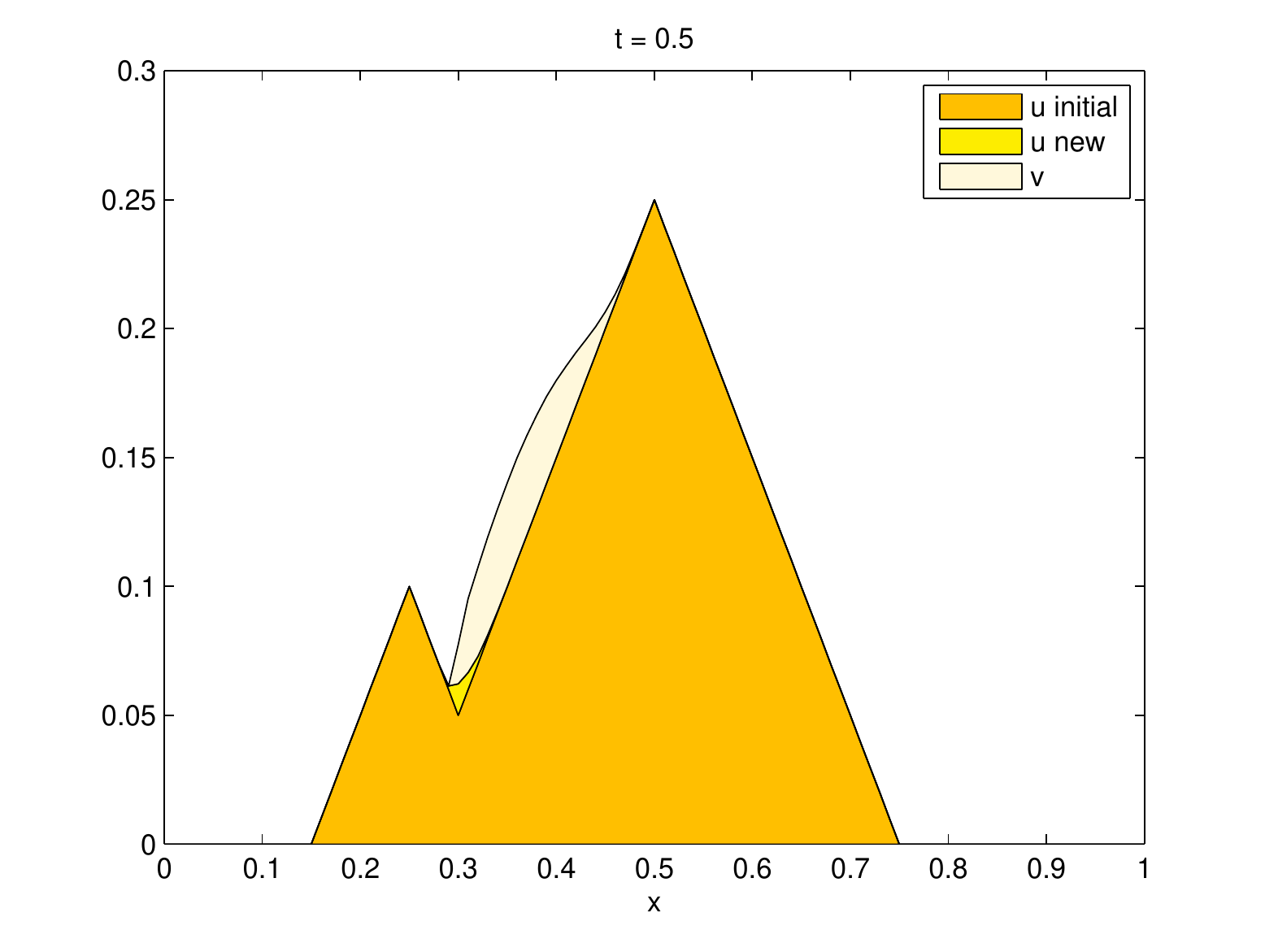}
\caption{Example 7.1: The case when the initial standing layer  has a cavity at $x=0.3$. Initial configuration (left) and configuration at t=0.5 (right).}
\label{cav1-2}
\end{figure}
\end{center}

\begin{center}
\begin{figure}
\includegraphics[width=7.5cm, height=4cm]{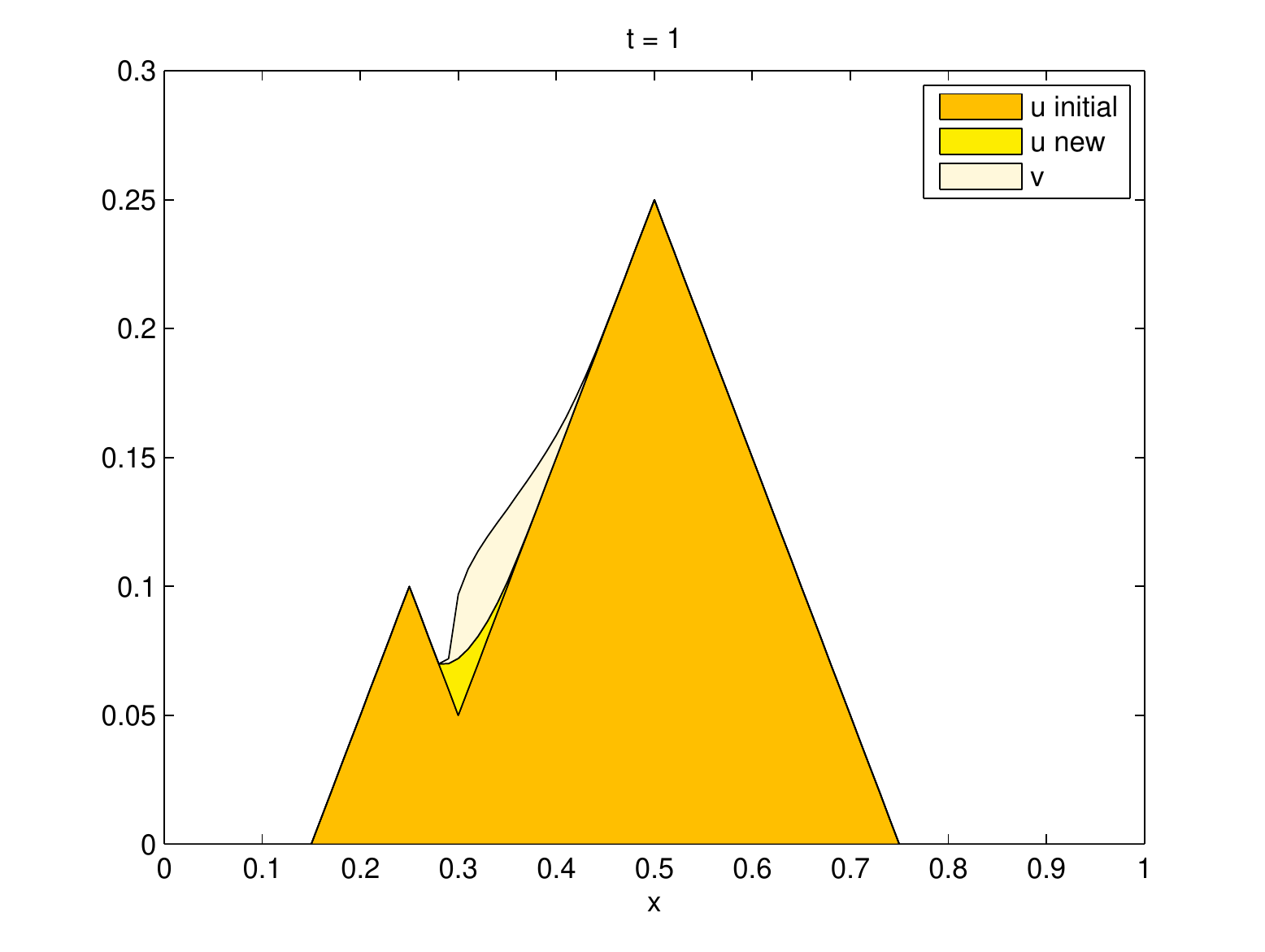}
\includegraphics[width=7.5cm, height=4cm]{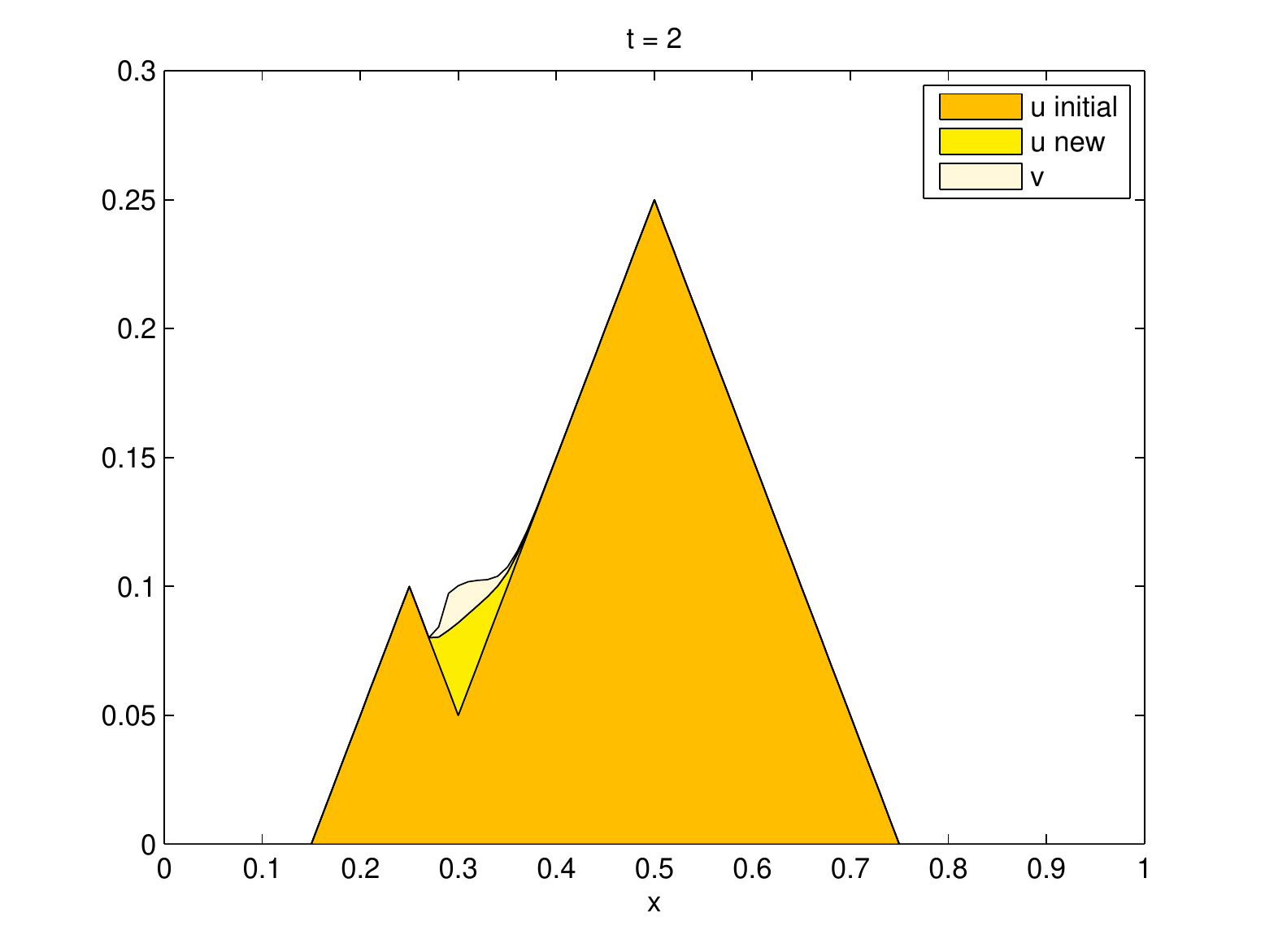}
\caption{Example 7.1: Configuration at t=1.0 (left) and t=2.0 (right).}
\label{cav3-4}
\end{figure}
\end{center}

\begin{center}
\begin{figure}
\includegraphics[width=7.5cm, height=4cm]{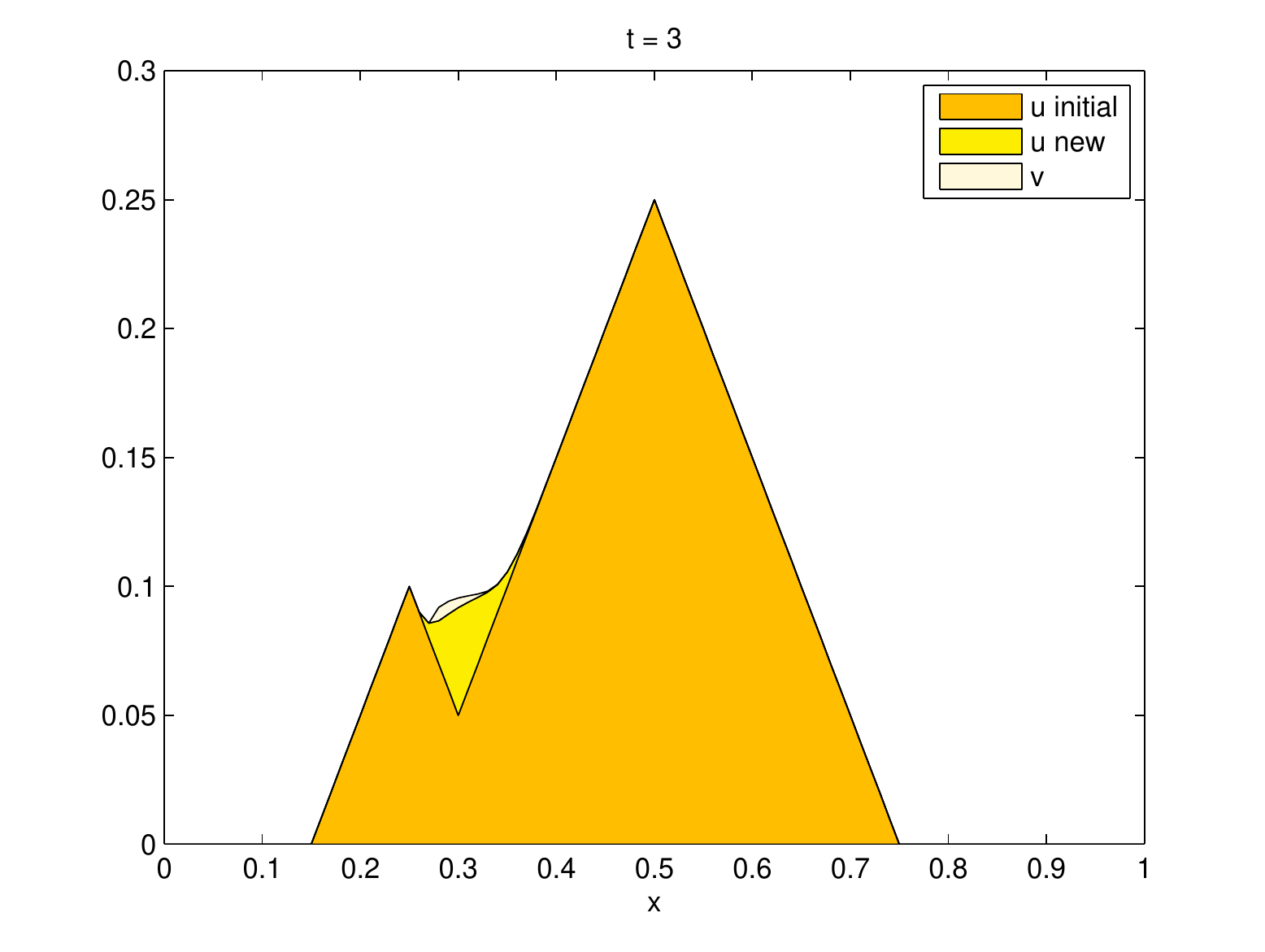}
\includegraphics[width=7.5cm, height=4cm]{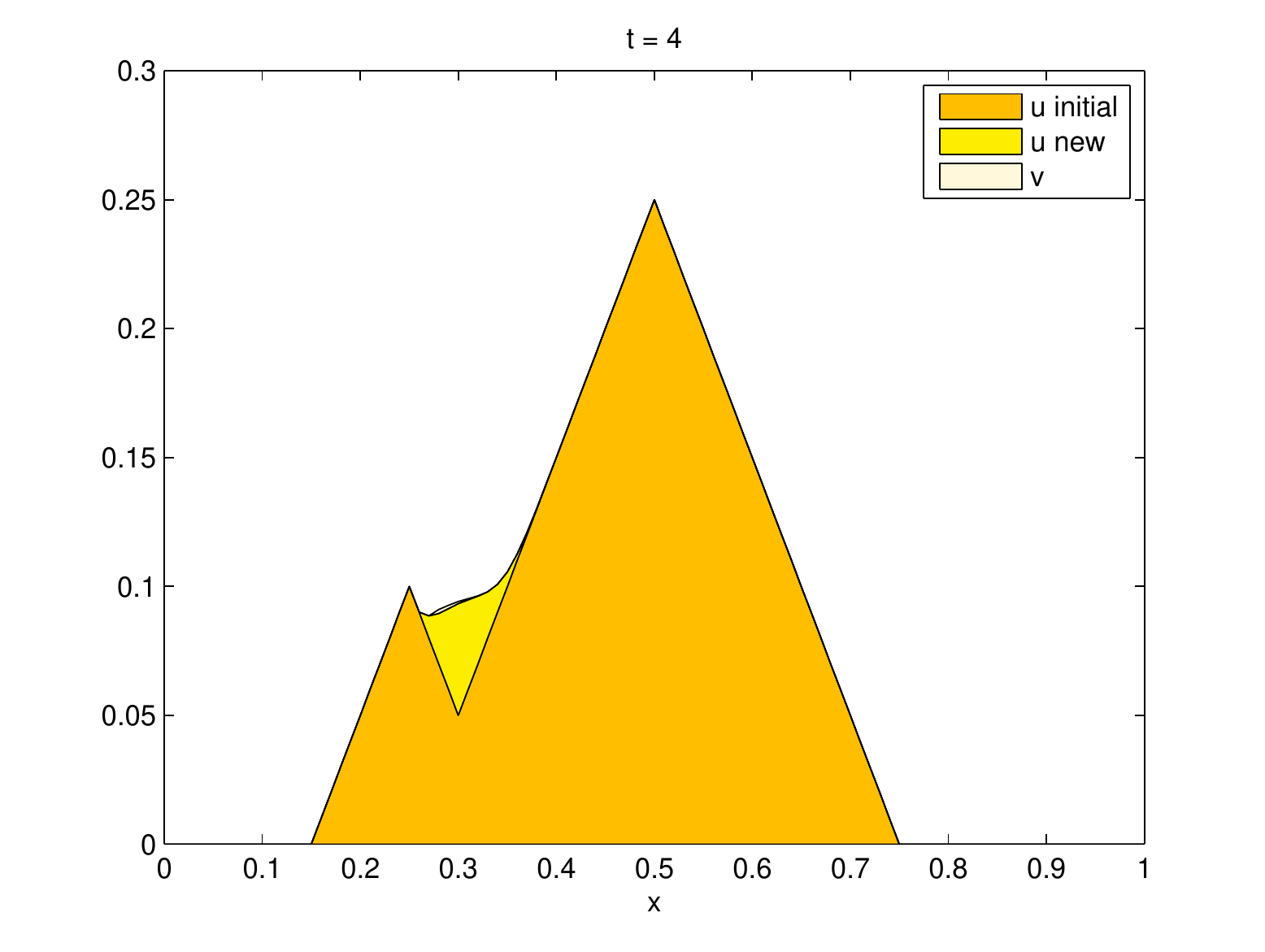}
\caption{Example 7.1: Configuration at t=3.0 (left) and t=4.0 (right).}
\label{cav5-6}
\end{figure}
\end{center}
\begin{center}
\begin{figure}
\includegraphics[width=7.5cm, height=4cm]{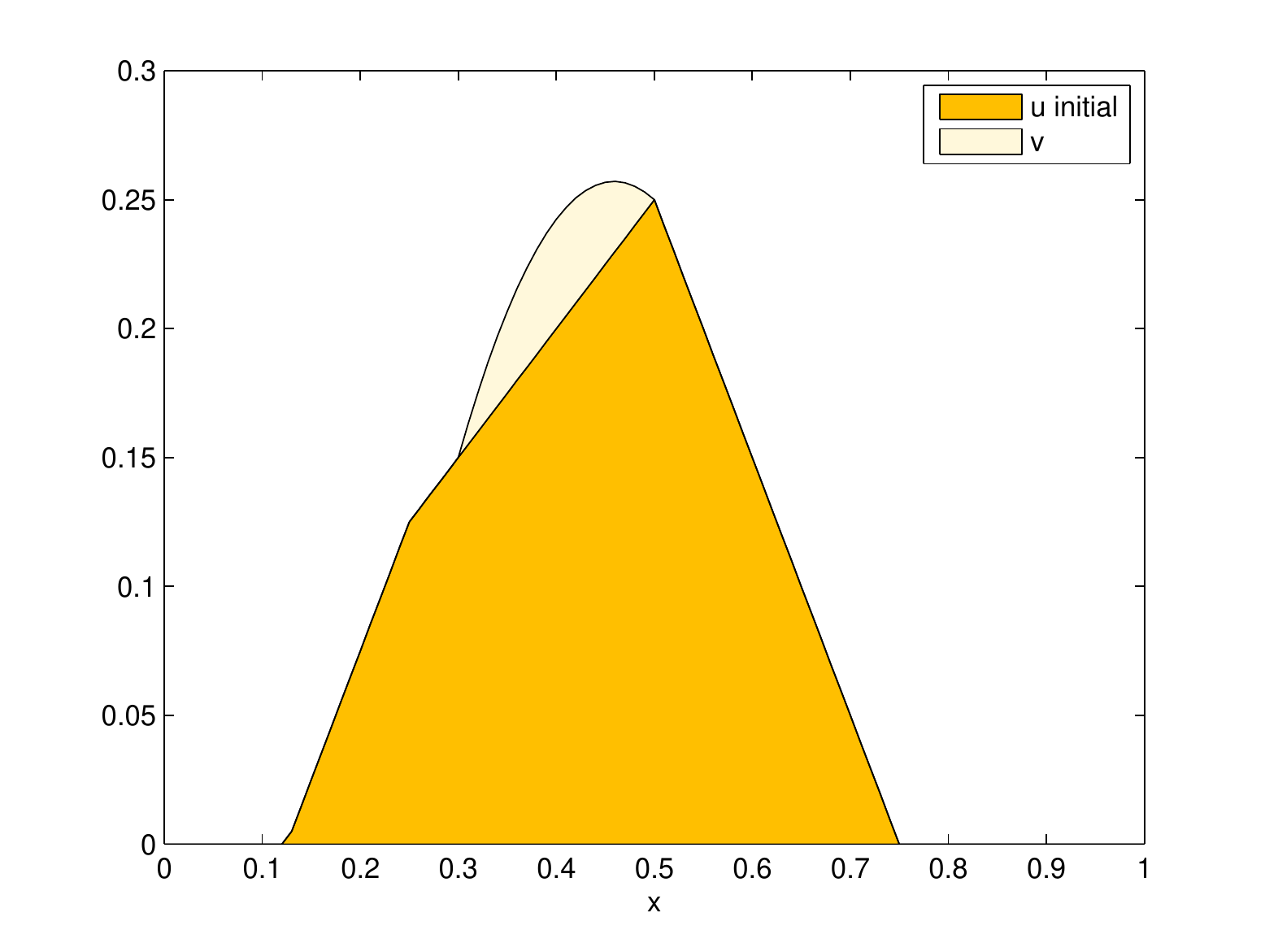}
\includegraphics[width=7.5cm, height=4cm]{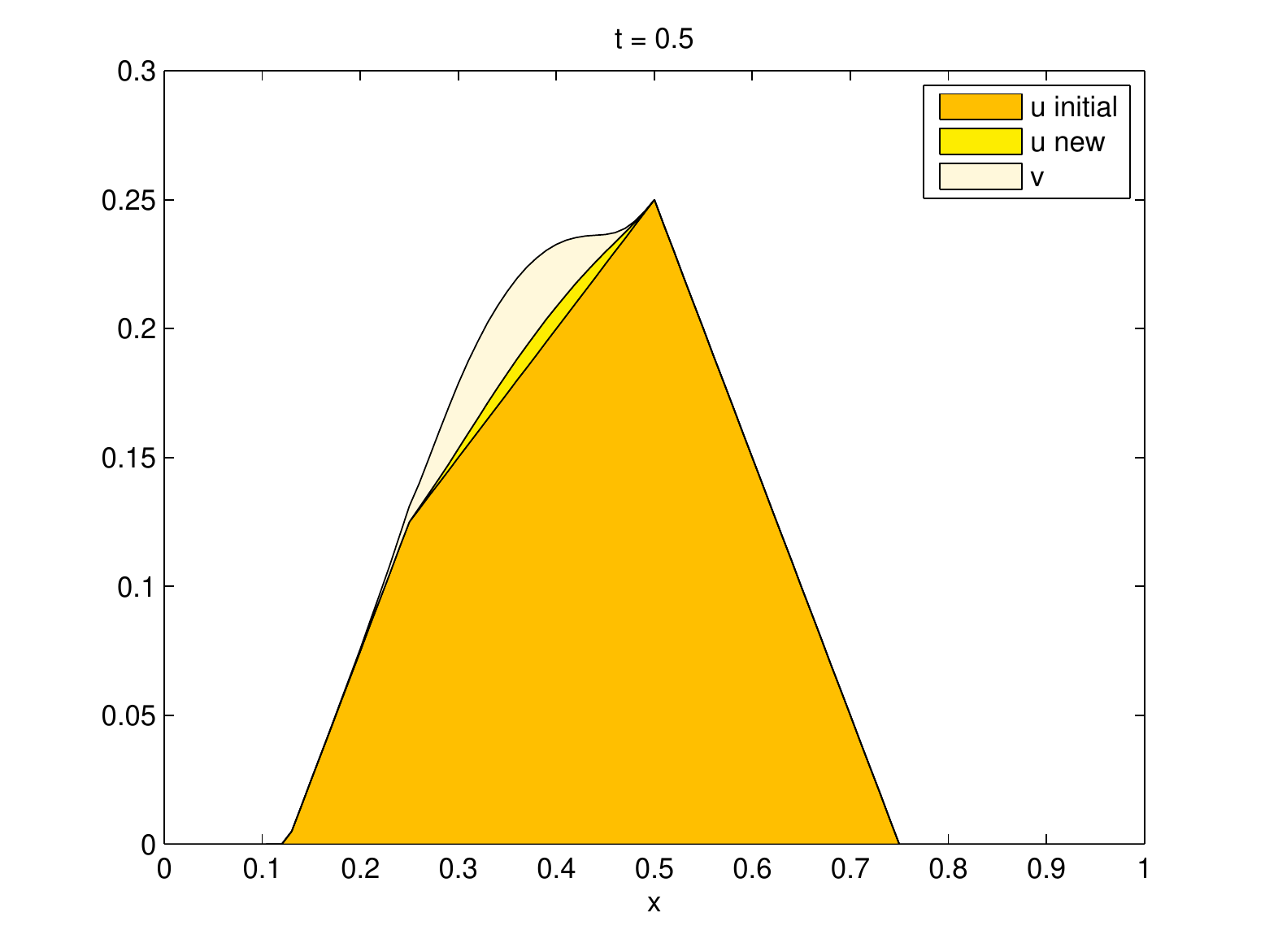}
\caption{Example 7.2 : The case when  the left slope of the initial standing layer changes from non-critical value to critical. Initial configuration (left) and t=0.5 (right).}
\label{cdown1-2}
\end{figure}
\end{center}
\begin{center}
\begin{figure}
\includegraphics[width=7.5cm, height=4cm]{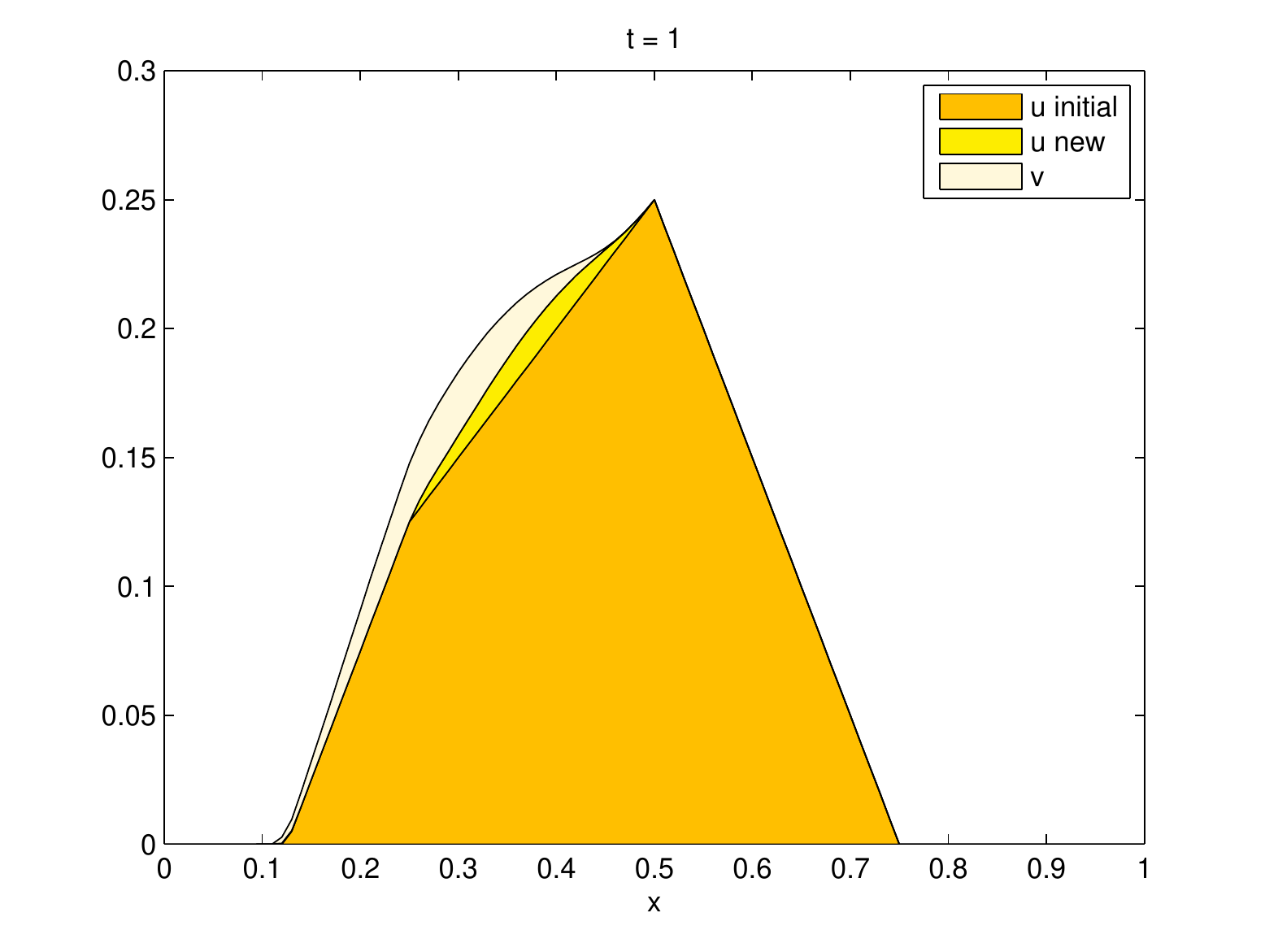}
\includegraphics[width=7.5cm, height=4cm]{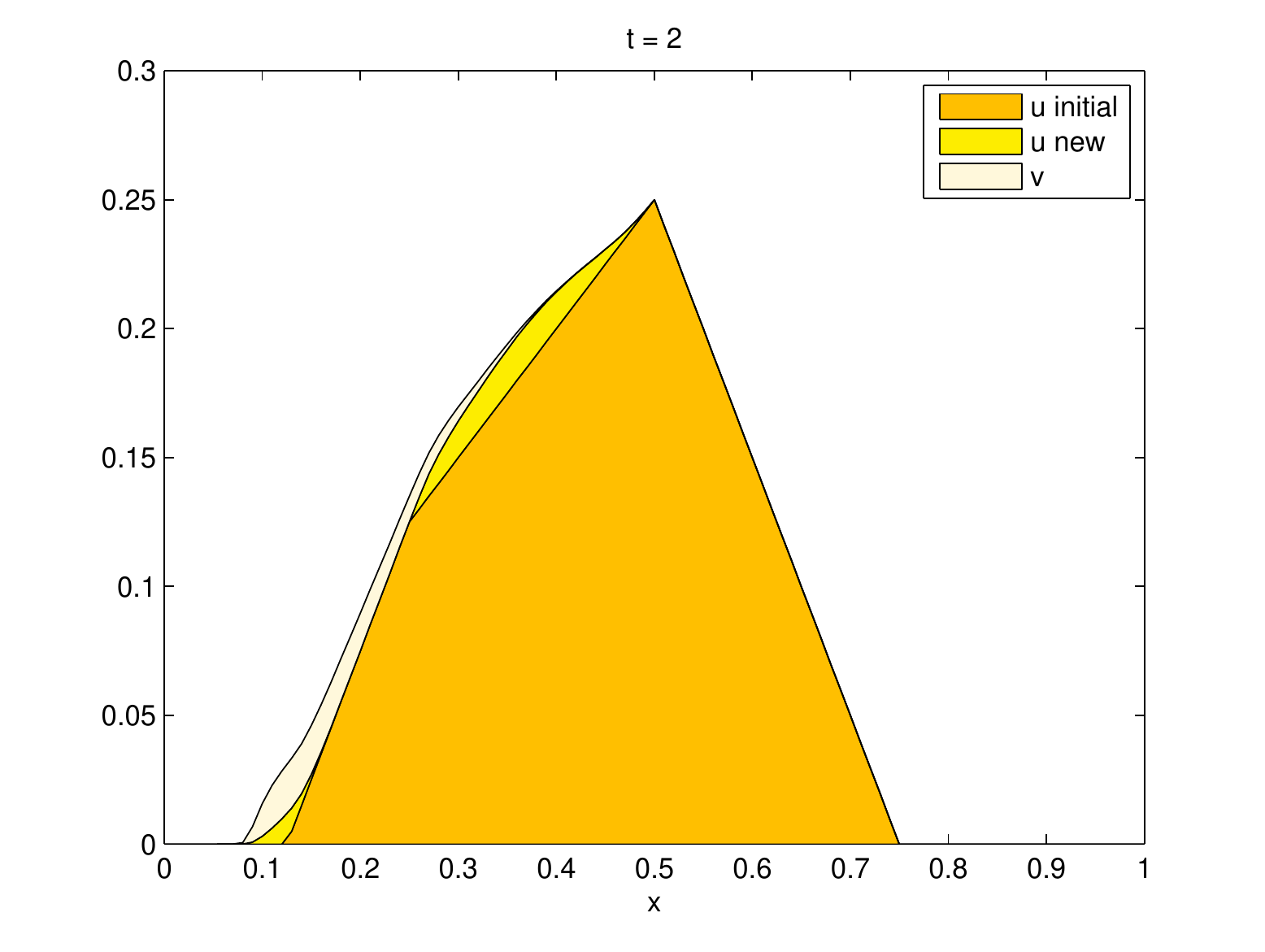}
\caption{Example 7.2: Configuration at t=1.0 (left) and t=2.0 (right).} 
\label{cdown3-4}
\end{figure}
\end{center}
\begin{center}
\begin{figure}
\includegraphics[width=7.5cm, height=4cm]{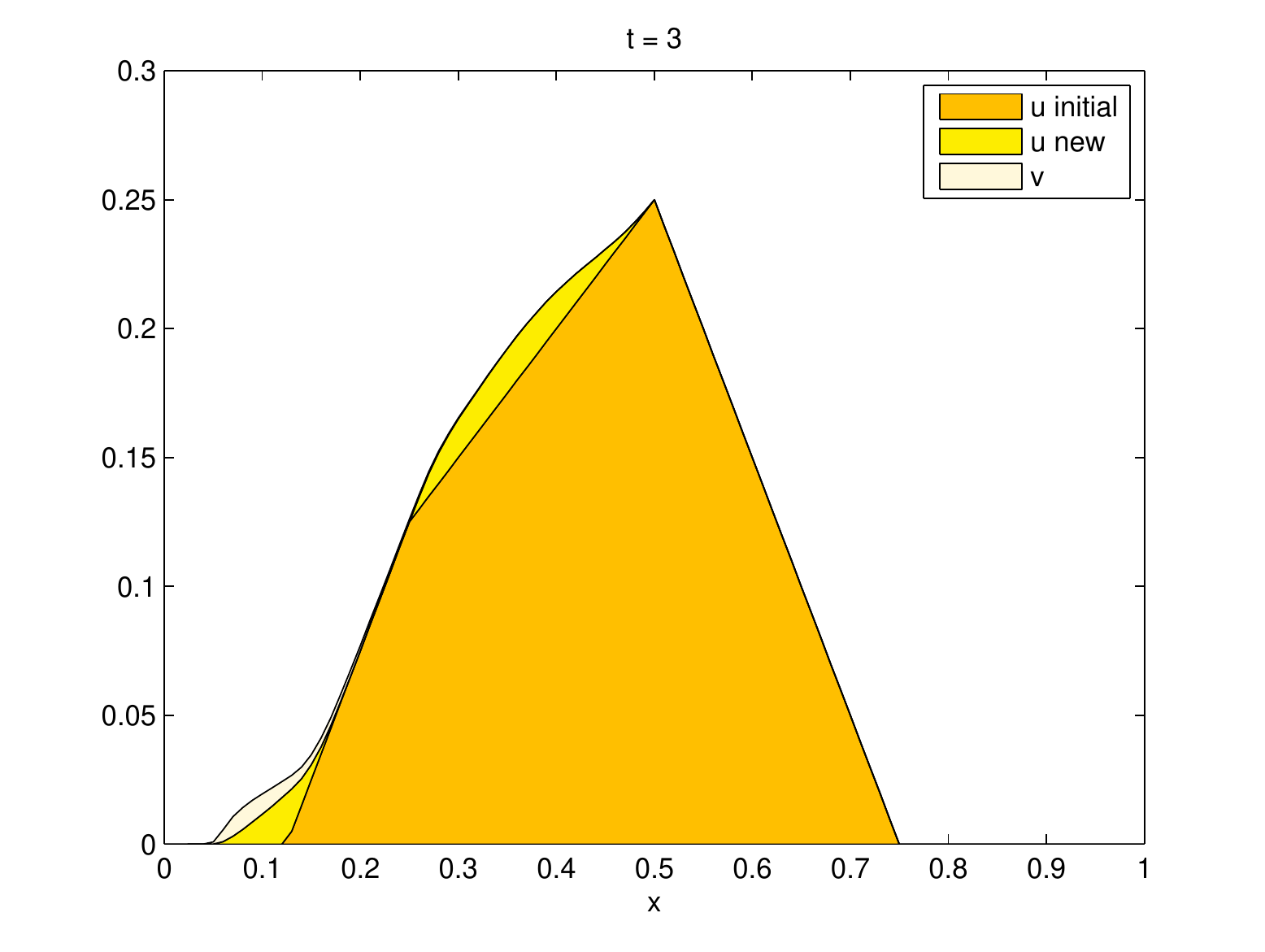}
\includegraphics[width=7.5cm, height=4cm]{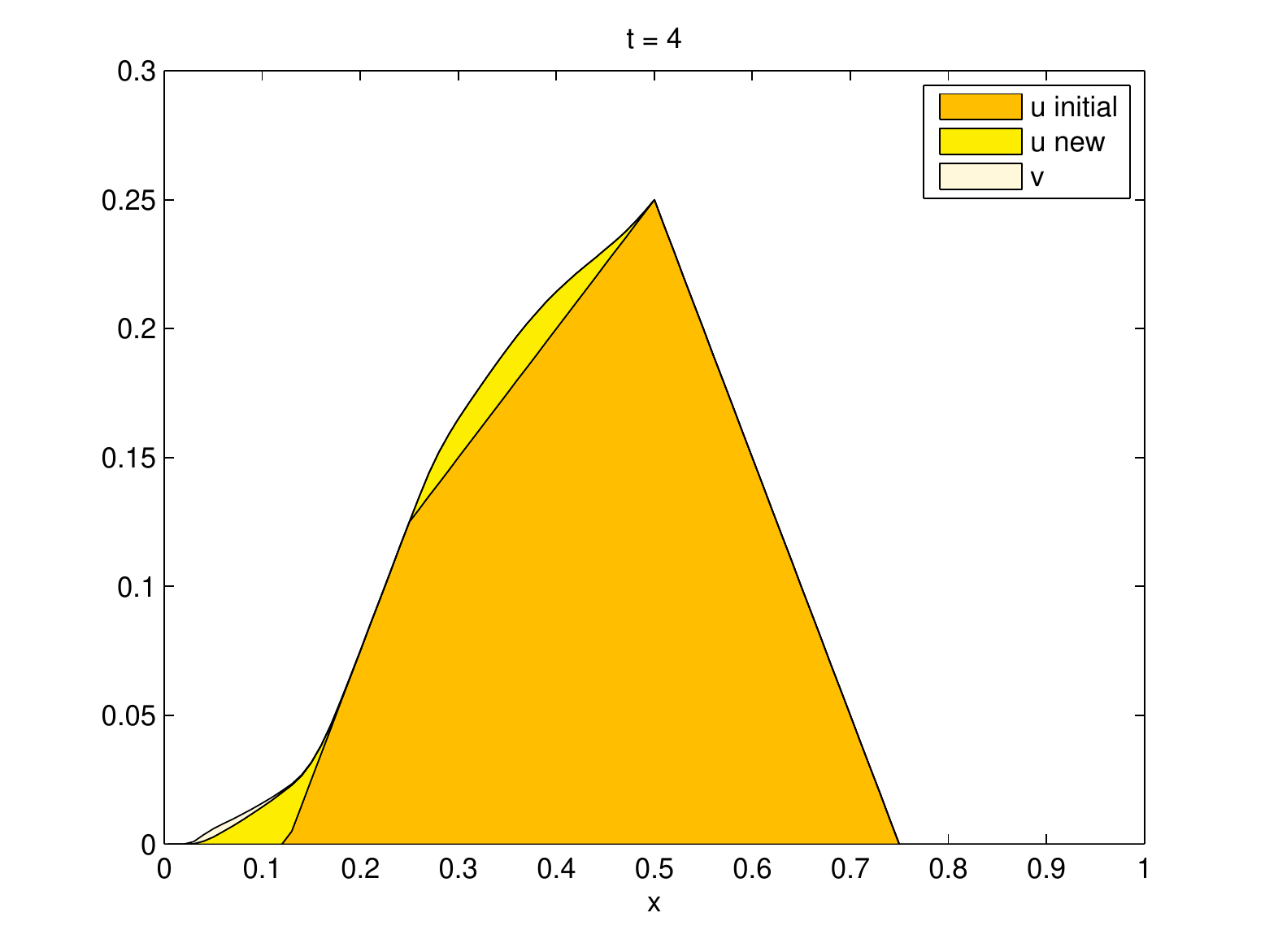}
\caption{Example 7.2: Configuration at t=3.0 (left) and t=4.0 (right).} 
\label{cdown5-6}
\end{figure}
\end{center}
 
 All the  models in the above (i.e., (1.1)-(1.2), (1.3)-(1.4)/(2.1)-(2.2), (2.4)-(2.5),(4.1)-(4.4)) do {\em not}  explicitly address the question of conservation of the total granular  matter. Our numerical results  indicate that, given the initial standing and rolling layers, this issue  is linked to the value of coefficient $ \beta$ in the equations for the rolling layer. 

 {\bf 8. 2-$D$ granular formation difference equations.} 
 We will now assume that $ u $ and $ v$ are functions of two spatial variables $ x$ and $y$: $ u = u (x,y, t), v = v (x,y, t)$.
 In the 1-$D$ model (4.1)-(4.4) we dealt with the flow of granular matter propagating along the $x$-axis only. In the 2-$D$ model we will  respectively deal with the flows propagating between the  nodes of the chosen mesh. We intend to use  the standard rectangular-shaped mesh on the $ xy$-plane and will propose a 2-$D$ discrete model which will describe the propagation of the granular matter between the respective adjacent nodes parallel to   the $ x$- and $y$-axes. As usual, one has certain freedom to chose the interpolation method for the spatial points which are not the nodes. 
  
\includegraphics[scale=0.31]{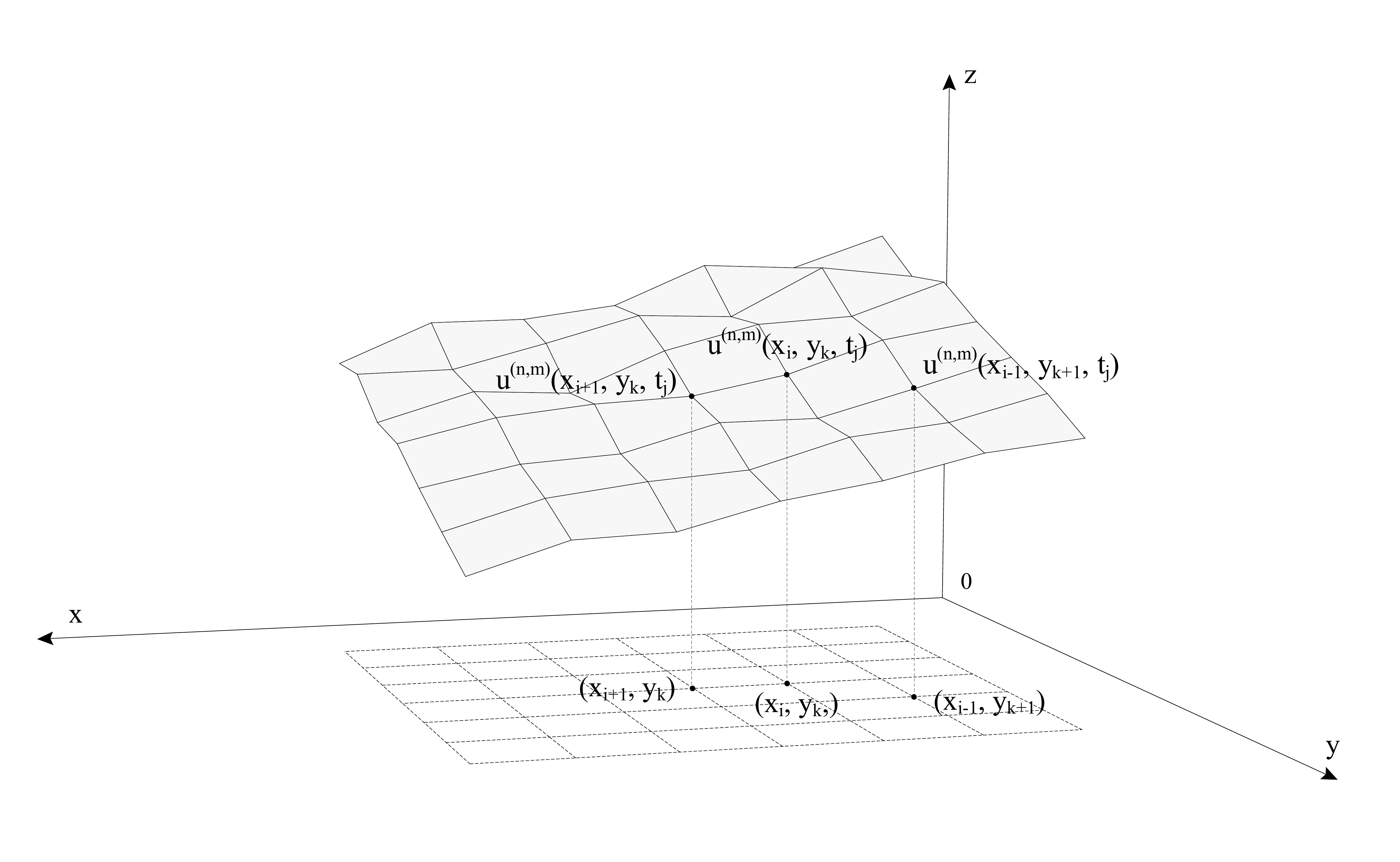} 
 
\begin{center}
Figure 15: Illustration of possible reconstruction of the standing layer in the 2-$D$ model.
\end{center}

We will preserve the notations of Section 4 wherever it will be possible. 
Let $ (x,y) \in [0, 1] \times [0,1] = \Omega$. Split  $ \Omega$  into $n \times n$ equal squares (one can  consider rectangles as well)   of size $\Delta_{n}  x  = \Delta_{n}  y $, and the time-interval $[0, T]$ into $m$ equal segments of size 
$ \Delta_m t$. The solution of the approximate system is represented by the collection of values denoted by $ \{ u^{(n, m)}_{i, k, j}, v^{(n,  m)}_{i, k, j}; i, k = 0, \ldots, n, j = 0, \ldots \}$. These values will define, depending on the chosen strategy for interpolation,  an   approximate solution in $\Omega$, denoted for the standing layer  by $  u^{(n,  m)} (x, y, t)$ and for the rolling layer by $  v^{(n,  m)} (x, y, t), (x,y)  \in [0,1] \times [0,1], t \in [0, T]$, see Fig 15 for illustration.

To approximate the time derivative at point $ (x_i, y_k, t_j)$, we use the same standard forward approximation as in Section 4. The same we will do for the $x$- and $y$-spatial derivative, e.g., 
$$
u_{x\pm}^{(n, m)} (x_i, y_k, t_j)  \; = \; \frac{u^{(n, m)} (x_{i\pm1}, y_k, t_{j}) - u^{(n, m)} (x_i, j_k, t_{j})}{\Delta_n x},   
$$
$$
u_{y\pm}^{(n, m)} (x_i, y_k, t_j)  \; = \; \frac{u^{(n, m)} (x_{i}, y_{k\pm1}, t_{j}) - u^{(n, m)} (x_i, j_{k}, t_{j})}{\Delta_n y}.   
$$

{\bf Difference equations for the  standing layer.} Its value will be determined by the {\em sum} of movements of the granular matter along the $x$- and $ y$-axes: 
$$
u^{(n, m)} (x_i, y_k, t_{j+1})   =  \; u_{(x)}^{(n, m)} (x_i, y_k, t_{j+1}) +\; u_{(y)}^{(n, m)} (x_i, y_k, t_{j+1}),
\eqno(8.1)$$
where $ u_{(x)}^{(n, m)} (x_i, y_k, t_{j+1})$ represents the propagation of this matter along the $x$-axis and is calculated for each $ y_k, k = 1, \ldots, n-1$  exactly as in (4.1) for the respectively modified notations in the list of variables: $ (x_i,t_j) $ in (4.1) should be replaced by $ (x_i, y_k, t_j)$. However, in this case each of  the  nodes $  (x_{i}, y_k) $  deals with  four directions,  i.e., along both the $x$- and $ y$-axes, instead of two as in the 1-$D$ case. In each term on the right of (8.1) (compare to (4.1))  we will need to use a respective  splitting coefficient.  
For example, if  at instant $ t_j$ at the node $ (x_i, y_k)$  in the negative direction along the $x$-axis we have the  steepest slope and the slopes in other three directions are strictly smaller,  then the respective splitting coefficient $ r_{i,k, j,x-}  = 1$. If it is one of two steepest slopes, $ r_{i,k, j,x-}  = 1/2 $, and so forth.

Denote 
$$
v_*^{(n, m)} (x_i, y_k,, t_{j})=  v^{(n, m)} (x_i, y_k, t_{j}) -\left(u^{(n, m)} (x_i, y_k, t_{j+1})  - u^{(n, m)} (x_i, y_k, t_{j}) \right), 
$$
where  the term in the large parenthesis describes the increase of the height of the standing layer during the time-interval $ (t_j, t_{j+1})$ due to the contribution from the rolling  layer available at time $ t_j$ (compare to  (4.1)-(4.4)). Thus,  $ v_*^{(n, m)} (x_i, y_k, t_{j}) $ describes the part of  the rolling layer at $x_i$ which will leave  this point rolling down a respective slope.  In other words, as in the 1-$D$ case, we assume that the contribution of the rolling layer to the standing layer on the interval $ [t_j, t_{j+1}]$ occurs at time $t_j$. Denote: 
$$
v_{*x\pm}^{(n, m)} (x_i, y_k, t_j)  \; = \; \frac{v_*^{(n, m)} (x_{i\pm1}, y_k, t_{j}) - v_*^{(n, m)} (x_i, t_{j})}{+\Delta_n x},   \;\;\;\; 
$$
$$
v_{*y\pm}^{(n, m)} (x_i, y_k, t_j)  \; = \; \frac{v_*^{(n, m)} (x_{i}, y_{k\pm1}, t_{j}) - v_*^{(n, m)} (x_i, y_k,t_{j})}{+\Delta_n x}.
$$

{\bf The difference equations for the  rolling layer.} In the above notations we obtain the following system of equations: 
$$
v^{(n, m)} (x_i,  y_k, t_{j+1})=  v_*^{(n, m)} (x_i,  y_k, t_{j}) \; +\; \Delta_m t \;   \beta \; ( \; {\mathbb F}_x \; + \;{\mathbb F}_y)
\eqno(8.2)$$
where  $ {\mathbb F}_x$  represents the motion of the rolling layer 
along the $x$-axis and $ {\mathbb F}_y$-- along the $y$-axis with the following ``correction''. Namely, the respective terms in (8.2) will have the splitting coefficients like  $r_{i,k, j,x\pm}$ in (8.1) in  the above.

{\bf The initial and boundary conditions} for equations (8.1)-(8.2)  are defined similar to those in  (4.1)-(4.2):
$$
  u^{(n, m)} (x_i, 0, t_j)  = u^{(n, m)} (x_k, 1, t_j) = u^{(n, m)} (0, y_k, t_j)  = u^{(n, m)} (1, y_k, t_j) = 0, \;\;\;\; j = 1, \ldots,
$$
$$
u^{(n, m)} (x_i, y_k, 0) = u_0 (x_i, y_k,) \; \geq 0, \;\; v^{(n, m)} (x_i,  y_k, 0) = v_0 (x_i,  y_k, 0) \; \geq 0, \;\; i, k = 0, \ldots, n,
\eqno(8.3)$$
$$
\max\{ \mid u_{x\pm}^{(n, m)} (x_i, y_k, 0) \mid, \; \mid u_{y\pm}^{(n, m)} (x_i, y_k, 0) \mid\}  \; \leq \;\alpha, \;\; \;\; i, k = 1, \ldots, n-1.
 \eqno(8.4)$$
 The convergence results of Section 6 can be extended to model (8.1)-(8.4) in a similar way.

 {\bf Connection to pde modeling.}  The difference equations (8.1)-(8.4) can be viewed as approximation of the following pde model    at the points $ (x, y, t)$  where $ u$ and $ v$ are continuously differentiable  and $ u_x \neq 0 \neq u_y$ (compare to (2.4)-(2.5)):
$$
v_t =   -u_t \; 
+\;  \beta \nabla u \cdot \nabla v, 
\eqno(8.5)$$
$$
u_t  \;  =  \;  \gamma  (\alpha  - \mid  \nabla u \mid), 
\eqno(8.6)$$
$$
u(x,0) = u_o (x) \geq 0, \;\;  v(x,0) = v_o (x) \geq 0, \;\; \max \{\mid u_{ox} \mid, \mid u_{oy} \mid \} \; \leq \alpha \;\;  \;\; x \in [0, 1],
$$
$$
u(x,t) \mid_{\partial \Omega} = 0, \;\;\;\; v(x,t) \mid_{\partial \Omega} = 0 \;\; t \in (0, T). 
$$
Note that our discrete model (8.1)-(8.2) is intended to restrict the maximum of directional derivatives of $u$ by the critical value along the directions parallel to the $x$- and $y$-axes (with further use of suitable interpolation elsewhere). When $ u$ is differentiable, this maximum is equal to $ \mid  \nabla u \mid$. Note that it is not the case when we deal with non-differentiable points of the standing layer such as the vertices  of cones. In this case, one can use the following equation in place of (8.6):
$$
u_t  \;  =  \;  \gamma  (\alpha  - \max_{s \in R^2} \{  \mid D_s u \mid), 
$$
where $ D_s u $ stands for the directional derivative of $ u$ in the direction $ s$.

{\bf Concluding remark on the discrete nature of models (4.1)-(4.4), (8.1)-(8.4).} One can view our discrete models (4.1)-(4.4) and (8.1)-(8.4) as ``unnecessarily complex''. Let us argue, however, that in order to calculate a solution of any highly nonlinear granular matter formation model, one will have to use {\em some} discrete model. Thus, the crux of the problem here  is the strategy on how to choose one, i.e.,  not just any abstract discretization, but a discretization which will preserve the physical properties of the process at hand. Models (4.1)-(4.4) and (8.1)-(8.4) are aimed exactly at this issue (see also an opening paragraph in Section 7 on the numerical aspects of various approximation  methods).

\end{document}